\documentclass[11pt]{article}
\usepackage{amsmath,amssymb}
\usepackage{pstricks}
\usepackage{pst-node}
\usepackage{pst-poly}
\usepackage{multido}
\usepackage{tikz}

\newtheorem{result}{Theorem}
\newtheorem{deduce}{Corollary}

\newtheorem{define}{Definition}
\newtheorem{support}{Lemma}
\newtheorem{propo}{Proposition}

\newtheorem{note}{Remark}

\newcommand{\qed}{%
\ifmmode % if math mode, assume display: omit penalty etc.
\else \leavevmode\unskip\penalty9999 \hbox{}\nobreak\hfill \fi
\quad\hbox{\qedsymbol}}
\newcommand{\openbox}{\leavevmode \hbox to.77778em{%
\hfil\vrule
\vbox to.675em{\hrule width.6em\vfil\hrule}%
\vrule\hfil}}
\newcommand{\qedsymbol}{\openbox}

\newcommand{\showgrid}{}
\newcommand{\gridon}{\renewcommand{\showgrid}{\psset{subgriddiv=1,griddots=10,gridlabels=6pt}\psgrid}}
\gridon

\begin{document}

\begin{center} {\bf\LARGE  Non-Crossing Perfect Matchings and Triangle-Free Geometric Graphs} \end{center}

\vskip8pt

\centerline{\large Hazim Michman  Trao$^{a}$, Gek L. Chia$^{b,e}$,  \  Niran Abbas Ali$^c$ \ and \   Adem Kilicman$^d$ }

\begin{center}
\itshape\small  $^{a, c, d}\/$Department of Mathematics, \\ Universiti Putra Malaysia, 43400 Serdang, Malaysia,  \\
 \vspace{1mm}
 $^{b}\/$Department of Mathematical and Actuarial Sciences, \\  Universiti Tunku Abdul Rahman, Sungai Long Campus,    Malaysia \\
\vspace{1mm}
 $^e\/$Institute of Mathematical Sciences, University of Malaya, \\ 50603 Kuala Lumpur,  Malaysia  \\
\end{center}

\begin{abstract}
 We study extremal type problem arising from the question:  What is the maximum number of edge-disjoint non-crossing perfect matchings on a set $S\/$ of $2n\/$ points in the plane such that their union is a triangle-free geometric graph?  We approach this problem by considering four different situations  of $S\/$.  In particular, in the general position, we obtain (i) a sufficient condition for the existence of $n\/$ edge-disjoint  non-crossing perfect matchings in the general position whose union is a maximal triangle-free geometric graph, and (ii) a lower bound  on the  number of  edge-disjoint non-crossing perfect matchings whose union is a  triangle free geometric graph.

 %We study extremal type problem arising from the question: What is the maximum number of edge disjoint non-crossing perfect matchings on a set $S\/$ of $2n\/$ points in the plane such that their union is a triangle-free geometric graph? We approach this problem by considering four different situations  of $S\/$, namely, when $S\/$ is in convex position, regular wheel position, general position, and a special case of general position (called $R\/$-position). In most of the situations, we also investigate the structures of the graphs obtained.
\end{abstract}

\vspace{1mm}
\section{Introduction}
Let $S$ be a set of points in general position in the plane (that is, no three points on a common line). By a {\em geometric graph\/} we mean a graph  $G=(S,E)$ where $S$ is a finite set of points ({\em vertices\/}) in general position in the plane and $E$ is a finite set of straight-line segments ({\em edges\/}) which joins vertices in $S$.\  In the case where all vertices in $S\/$ are  in a convex position, the  graph is called a {\em convex geometric graph\/}.

\vspace{1mm} Two subgraphs in a geometric graph are said to be {\em edge-disjoint\/} if they have have no edge in common.  Two edges in a geometric graph are said to {\em  cross\/}  each other if they have a common  point which is interior to both of them.  A subgraph of a geometric  is said to be {\em non-crossing\/} (or {\em plane\/})  if its edges do not cross each other.

\vspace{1mm} A geometric graph $G\/$ is said  to be {\em triangle-free\/} if it contains no complete subgraph with three vertices.

\vspace{1mm} By a {\em matching\/} in a graph $G\/$ we mean a subgraph of $G\/$ in which every vertex has degree at most one; and a  matching is {\em perfect\/} if every vertex has degree equal to one.

\vspace{1mm}
 Matchings in graphs have received a  considerable amount of attention in graph theory (see \cite{lp:refer}).  Meanwhile matching problems have also been considered in the context of triangle-free graphs (see \cite{bhhsv:refer} and \cite{lm:refer}).

%\vspace{1mm} It is well-known that the maximum matchings problem and triangle-free graphs play an important roles in graph theory and combinatorial optimization since they have a wide range of applications. The aspects: counting non-crossing perfect matchings in graphs, and the structure of triangle-free graphs have been merged in many works (see \cite{bhhsv:refer} and \cite{lm:refer} ). In particular, in \cite{bhhsv:refer} the authors proved that every n-vertex triangle-free graph has at most $3^{n/3} \approx 1.4423^{n}\/$ maximal induced matchings.}

\vspace{1mm}
For the case of geometric graphs, matchings have also received much attention in the literature (see \cite{abz:refer}, \cite{aic:refer}, \cite{ata:refer}, \cite{bbms:refer} and \cite{hhn:refer}  for example). In particular, in \cite{bbms:refer}, the authors consider the problem of packing a  maximum number of  edge-disjoint non-crossing perfect matchings into a convex complete graph of even order.   In turns out, as we shall show in Section \ref{convexposition}  that   when  the maximum number of edge-disjoint non-crossing perfect matchings is packed into a set $S\/$ (of even number of points in convex position),   the result is a maximal triangle-free graph   (see Theorem \ref{c2n-n}).

\vspace{1mm}
For the case where $S\/$ is a set of $2n\/$ points in regular wheel configuration in the plane (Section \ref{wheel}), we show that  there exist at most  $\lceil\frac{(b-1)n}{b}\rceil -(b-2)$  edge-disjoint non-crossing perfect matchings whose union yields a triangle-free geometric graph. Here $b\/$ denotes the number of boundary edges in each such perfect matching, and $n \geq b^2-1\/$ (see Theorem \ref{b=2,3}). It is noted that the bound is best possible when $b=3\/$. Also, when $b=2\/$, we show that the resulting graph obtained is unique (Theorem \ref{w-oddn2}).

\vspace{1mm}
For the case $S\/$ is a set of $2n\/$ points  in  general position in the plane where $n=2^k+h$ (and $0\leq h <2^k$),  we show that there exist at least $k$ edge-disjoint non-crossing perfect matchings whose union is a triangle-free geometric graph. The proof makes use of an algorithm (Algorithm $(A)\/$) in Section \ref{general} which is modified from an algorithm given in \cite{ata:refer}.

\vspace{1mm}  In Section  \ref{rposition}, we present a  sufficient condition for the existence of $n\/$ edge-disjoint non-crossing perfect matchings (on a given set of $2n\/$ points in the general position) whose union is a maximal triangle-free geometric graph (Theorem \ref{r-post}).

\vspace{1mm}
The case where the resulting graph (that is, the triangle-free geometric graph obtained from the union of $k\/$ non-crossing perfect matchings) is a geometric plane graph is treated in the final section.

\vspace{1mm} The following result will be used quite often in the rest of the paper. Let $CH(S)\/$ denote the convex hull of a set $S\/$.

\begin{propo} \label{s-even}
Let $S\/$ be a set of $2n\/$ points, and $F\/$ be a non-crossing perfect matching on $S\/$ where $n \geq 2\/$. Then each edge $uv\/$ in $F\/$ such that $u, v \in CH(S)\/$ divides $S - \{u, v\}\/$ into two parts each having an even number of vertices.
\end{propo}

\vspace{1mm}  \noindent
{\bf Proof:} Suppose there is an edge $uv\/$ in $F\/$ where $u, v \in CH(S)\/$ that divides $S - \{u, v\}\/$ into two parts each having an odd number of vertices (since $S\/$ has $2n\/$ points). But then $F\/$ cannot match  all the vertices of each part of $S\/$ without crossing the edge $uv\/$ (since each part has an odd number of vertices). This contradiction proves the lemma.   \qed

\section{Points in convex position} \label{convexposition}

\vspace{1mm} In \cite{bbms:refer}, the authors prove that for  a set $S\/$ of $2n\/$ points in  convex position, there are at most $n\/$  edge-disjoint non-crossing  perfect matchings on $S\/$.   We shall show that when $n\/$ such   edge-disjoint non-crossing perfect matchings are packed into $S\/$, the resulting convex geometric graph is the unique maximal triangle-free graph $C_{2n,n}\/$  defined below.

\vspace{1mm}
Let $G\/$ be a convex graph whose vertices are arranged in cyclic order $u_0, u_1, \ldots, u_{2n-1}\/$. Let $u_iu_{i+1}\/$ be an edge in $G\/$. An edge $u_ku_l\/$ is said to be {\em $p\/$-parallel\/} to $u_iu_{i+1}\/$ if  $k+l \equiv 2i+1 \ (mod \ 2n)\/$. Two edges  are said to be {\em $p\/$-parallel\/} if they are both $p\/$-parallel to the same boundary edge.

\vspace{1mm}   Let the  vertices of $C_{2n,n}\/$ be denoted by  $v_0, v_1, v_2, \ldots, v_{2n-1}\/$ arranged in cyclic order. For each $i = 0, 1, \ldots, n-1\/$ let $F_i\/$ denote the set of all edges $p\/$-parallel to the edge $v_0v_{2i+1}\/$.  Then clearly each $F_i\/$ is a non-crossing perfect matching on $S= \{v_0, v_1, v_2, \ldots, v_{2n-1}\}\/$. Moreover $E(F_i) \cap E(F_j) = \emptyset\/$ for $i \neq j\/$. We may color the vertices $v_0, v_2, \ldots, v_{2n-2}\/$ with one color and the vertices $v_1, v_3, \ldots, v_{2n-1}\/$ with another color.  Let  $ E(C_{2n,n}) = \bigcup_{i=0} ^{n-1} E(F_i)\/$. It is easy to see that  $C_{2n,n}\/$ is a triangle-free graph  having $2n\/$ vertices and $n^2\/$ edges (since no edge joins two vertices of the same color). By Turan's theorem, $C_{2n,n}\/$ is a maximal triangle-free graph.   Figure \ref{c126} depicts the graph $C_{12,6}\/$.

\begin{figure}[htb]
\centering
\begin{tikzpicture}
       \coordinate (v1) at (3,0);\node at (v1) {\textbullet};\node[above right] at (v1) {${v_1}$};
       \coordinate (v2) at (2.6,1.5);\node at (v2) {\textbullet};\node[above right] at (v2) {${v_2}$};
       \coordinate (v3) at (1.5,2.6);\node at (v3) {\textbullet};\node[above right] at (v3) {${v_3}$};
       \coordinate (v4) at (0,3);\node at (v4) {\textbullet};\node[above] at (v4) {${v_4}$};
       \coordinate (v5) at (-1.5,2.6);\node at (v5) {\textbullet};\node[above left] at (v5) {${v_5}$};
       \coordinate (v6) at (-2.6,1.5);\node at (v6) {\textbullet};\node[above left] at (v6) {${v_6}$};
       \coordinate (v7) at (-3,0);\node at (v7) {\textbullet};\node[above left] at (v7) {${v_7}$};
       \coordinate (v8) at (-2.6,-1.5);\node at (v8) {\textbullet};\node[below left] at (v8) {${v_8}$};
       \coordinate (v9) at (-1.5,-2.6);\node at (v9) {\textbullet};\node[below left] at (v9) {${v_9}$};
       \coordinate (v10) at (0,-3);\node at (v10) {\textbullet};\node[below] at (v10) {${v_{10}}$};
       \coordinate (v11) at (1.5,-2.6);\node at (v11) {\textbullet};\node[below] at (v11) {${v_{11}}$};
       \coordinate (v0) at (2.6,-1.5);\node at (v0) {\textbullet};\node[below right] at (v0) {${v_{0}}$};
         %lines
       \draw [line width=1,red](v1) -- (v2);\draw [line width=1,red](v0) -- (v3);\draw [line width=1,red](v11) -- (v4);
       \draw [line width=1,red](v10) -- (v5);\draw [line width=1,red](v9) -- (v6); \draw [line width=1,red](v8) -- (v7);
       \draw [line width=1,blue](v2) -- (v3); \draw [line width=1,blue](v1) -- (v4); \draw [line width=1,blue](v0) -- (v5);
       \draw [line width=1,blue](v11) -- (v6);\draw [line width=1,blue](v10) -- (v7);\draw [line width=1,blue](v9) -- (v8);
       \draw [line width=1,green](v3) -- (v4); \draw [line width=1,green](v2) -- (v5); \draw [line width=1,green](v1) -- (v6);
       \draw [line width=1,green](v0) -- (v7); \draw [line width=1,green](v11) -- (v8); \draw [line width=1,green](v10) -- (v9);
       \draw [line width=1,orange](v4) -- (v5); \draw [line width=1,orange](v3) -- (v6); \draw [line width=1,orange](v2) -- (v7);
       \draw [line width=1,orange](v1) -- (v8); \draw [line width=1,orange](v0) -- (v9); \draw [line width=1,orange](v11) -- (v10);
       \draw [line width=1,yellow](v5) -- (v6); \draw [line width=1,yellow](v4) -- (v7); \draw [line width=1,yellow](v3) -- (v8);
       \draw [line width=1,yellow](v2) -- (v9); \draw [line width=1,yellow](v1) -- (v10); \draw [line width=1,yellow](v0) -- (v11);
       \draw [line width=1,purple](v6) -- (v7); \draw [line width=1,purple](v5) -- (v8); \draw [line width=1,purple](v4) -- (v9);
       \draw [line width=1,purple](v3) -- (v10); \draw [line width=1,purple](v2) -- (v11); \draw [line width=1,purple](v1) -- (v0);

      \node at (v1) {\textbullet};\node at (v2) {\textbullet}; \node at (v3) {\textbullet}; \node at (v4) {\textbullet};
      \node at (v5) {\textbullet}; \node at (v6) {\textbullet}; \node at (v7) {\textbullet};\node at (v8) {\textbullet};
      \node at (v9) {\textbullet};\node at (v10) {\textbullet}; \node at (v11) {\textbullet};\node at (v0) {\textbullet};

\end{tikzpicture}
\caption{$C_{12,6}\/$}  \label{c126}
\end{figure}

\vspace{1mm}
The following result shows that any  $n\/$  edge-disjoint non-crossing  perfect matchings on a set of $2n\/$ points in convex position  yields the same  maximal triangle-free graph  $C_{2n,n}\/$.

\vspace{1mm} Suppose $S= \{v_0, v_1, v_2, \ldots, v_{m-1}\}\/$ is a set of $m\/$ points in the convex position. Edges of the form $v_iv_{i+1}\/$, $i=0,1,2,  \ldots, m-1\/$ are called the {\em boundary edges\/} of $S\/$.

\vspace{1mm}
\begin{result} \label{c2n-n}
Let $S\/$ be a set of $2n\/$ points in the convex position on the plane where $n \geq 2\/$. Suppose $ F_1, F_2, \ldots, F_n\/$ are $n\/$  edge-disjoint non-crossing  perfect matchings  on $S\/$. Then $\bigcup_{i=1} ^n F_i\/$   is   $C_{2n,n}\/$.
\end{result}

\vspace{1mm}  \noindent
{\bf Proof:} Because every  non-crossing perfect matching on $S\/$ has at least $2\/$ boundary edges on $S\/$ (see \cite{hhn:refer}), it follows that  $F_i\/$ contains exactly $2\/$ boundary edges on $S\/$ for each $i =1, 2, \ldots, n\/$.

\vspace{1mm}
Take any boundary edge on $S\/$, say $v_0v_1\/$. Then $v_0v_1\/$ belongs to one of the given  perfect matchings, say $v_0v_1 \in E(F_1)\/$. We assert that all edges in $F_1\/$ are $p\/$-parallel to $v_0v_1\/$.

\vspace{1mm} Suppose on the contrary that $v_iv_j \in E(F_1)\/$ is such that $v_0v_1\/$ and $v_iv_j\/$ are not $p\/$-parallel for some  $i < j\/$.  Then we have $|\{ v_2, \ldots , v_{i-1}\}| \neq |\{ v_{j+1}, \ldots, v_{2n-1}\}|\/$. If $j-i\/$ is an even integer, then there is an edge say $v_rv_s\/$ in $F_1\/$ with  $r \in \{i+1, i+2, \ldots, j-1\}\/$ and $s \not \in \{i, i+1, i+2, \ldots,  j\}\/$ implying $v_iv_j, v_rv_s\/$ are crossing edges in $F_1\/$.

\vspace{1mm} Hence $j-i\/$ is an odd integer and this implies that there is a boundary edge $v_tv_{t+1}\/$ on $S\/$ where $i <t < j\/$.

\vspace{1mm}
 Assume without loss of generality that $|\{ v_2, \ldots , v_{i-1}\}| < |\{ v_{j+1}, \ldots, v_{2n-1}\}|\/$. Since $F_1\/$ is a perfect matching on $S\/$, this implies that $F_1\/$ would have another boundary edge $v_kv_{k+1}\/$  on $S\/$ with $v_k, v_{k+1} \in   \{ v_{j+1}, \ldots, v_{2n-1}\}\/$, a contradiction.
 %But this is a contradiction since $F_1\/$ would have another boundary edge $v_kv_{k+1}\/$  on $S\/$ where $1 < k < i\/$ or $j < k < 2n-1\/$ (as $F_1\/$ contains no crossing edges).

\vspace{1mm}
Clearly, there exist $n\/$ boundary edges on $S\/$ where any two are non $p\/$-parallel. We can take these $n\/$ boundary edges to be $v_0v_1, v_1v_2, \ldots, v_{n-1}v_n\/$ and  assume that $v_{i-1}v_i  \in E(F_i)\/$, $i=1, 2, \ldots, n\/$. By the preceding argument, each $F_i\/$ yields a set of $n\/$ non-crossing edges. Since $v_0v_{2i+1} \in E(F_{i+1})\/$, $i=0, 1, \ldots, n-1\/$, we conclude that  $\bigcup_{i=0} ^{n-1} F_i\/$   is the maximal triangle-free convex geometric graph $C_{2n,n}\/$.     \qed

\vspace{2mm}
\begin{note}
Note that the bound on $n\/$ in Theorem \ref{c2n-n} is tight. Note also that it is possible  to construct  a maximal triangle-free geometric graph $G\/$ with interior points and $G\/$ is isomorphic to $C_{2n,n}\/$. An example of such geometric graph is shown in Figure \ref{sufficint}.
\end{note}

\vspace{1mm}
\begin{propo}  \label{separate-k}
Suppose $G\/$ is a convex geometric graph on $2n\/$ vertices and $G\/$ is the union of $n\/$ edge-disjoint non-crossing perfect matchings where $n \geq 3\/$. Let $K\/$ be a set of consecutive boundary vertices in $CH(G)\/$.  Then there exists a diagonal edge of $G\/$ that separates $K\/$ from the rest of the vertices of $G\/$ if and only if $|K|\/$ is even.
\end{propo}

\vspace{1mm}  \noindent
{\bf Proof:}  The necessity follows from Proposition \ref{s-even}.

\vspace{1mm} To show the sufficiency, we may assume without loss of generality that $ K=\{ v_1, v_2, \ldots, \linebreak v_k\}\/$ (since we can relabel the vertices by making a cyclic shift on the indices) and $k\/$ is even. Suppose there is no diagonal edge that separates $K\/$. Then $v_0v_{k+1}\/$ is not an edge of $G\/$ which means that $v_0, v_{k+1}\/$ belong to the same partite set of $G\/$ and hence $k+1\/$ is even. But this contradicts the fact that $|K| =k\/$ is even.   \qed

%newProposition

\section{Points in regular wheel configuration} \label{wheel}

A set $S\/$ of $m\/$ points is said to be in {\em regular wheel configuration\/} if $m-1\/$ of its points are regularly spaced on a circle $C\/$ with one point $x\/$ in the center of $C\/$. We call $x\/$ the {\em center} of $S\/$. Note that those vertices in $C\/$  are the convex hull of S.  An edge of the form $xv\/$ is called a {\em radial edge\/}; all other edges are called {\em non-radial edges\/}. Note that in this case, every perfect matching on $S\/$ contains a radial edge. By a {\em radial vertex\/} we mean a vertex that is incident with a radial edge.

\vspace{1mm} The following lemma will be used in the proofs of Theorems \ref{b=2,3} and \ref{w-oddn2} regarding the constructions of perfect matchings.

\vspace{1mm}
 Suppose $m =2n\/$ and $n\/$ is odd. Let $u_iu_{i+1}\/$ be an edge in $CH(G)\/$. (i) If $i \in \{(n-1)/2, (n+1)/2, \ldots, n-1\}\/$, then an edge $u_ku_l\/$ is said to be {\em $p_1\/$-parallel\/} to $u_iu_{i+1}\/$ if  $k+l \equiv 2i+1 \ (mod \ 2n)\/$ where $k+l \leq 2n-1\/$.  (ii) If  $i \in \{3(n-1)/2, (3n-1)/2, \ldots, 2n-2\}\/$, then an edge $u_ku_l\/$ is said to be {\em $p_2\/$-parallel\/} to $u_iu_{i+1}\/$ if  $k+l \equiv 2i+1 \ (mod \ 2n-1)\/$ where $k+l \leq 2n-2\/$.

\begin{support} \label{wheel-l1}
Let $S\/$ be a set of $2n\/$ points in regular wheel configuration  in the plane where $n \geq 3\/$ is odd. Suppose $F\/$ is a non-crossing perfect matching on $S\/$.

(i) Then $F\/$ has at least two boundary edges.

(ii) In the case that $F\/$ has only two boundary edges $e_1, e_2\/$, every non-radial edge in $F\/$ is $p_i\/$-parallel to either $e_1\/$ or $e_2\/$ for some $i \in \{1, 2\}\/$.

%(iii) {\bf In the case that $F\/$ has only three boundary edges $e_1, e_2\/$ and $e_3\/$, any  non-radial edge in $F\/$ is $p\/$-parallel to $e_i\/$ for some $i \in \{1, 2, 3\}\/$ ???  }
\end{support}

\vspace{1mm}  \noindent
{\bf Proof:} Suppose $x\/$ is the center of $S\/$ and $v_0, v_1,  \cdots, v_{2n-2} \/$ are the vertices of the circle $C\/$. Assume that $xv_0\/$ is the radial edge of $F\/$. Then $S - \{x, v_0\}\/$ is a set of $2n-2\/$ points in convex position. This means that $F -v_0x\/$ is a non-crossing perfect matching in $S - \{x, v_0\}\/$ and hence it has at least two boundary edges. This proves (i).

\vspace{1mm}
If $n = 3\/$, then (ii) is clearly true. Hence assume that $n \geq 5\/$.

\vspace{1mm} Clearly, $v_1v_i \in E(F)\/$ for some $i \in \{2, 4, \ldots, 2\lfloor n/2 \rfloor\}\/$.  But this means that $\{v_1, v_2, \ldots, \linebreak v_i\}\/$ contains a non-crossing perfect matchings $F_1\/$ with a boundary edge $e_1\/$. Moreover any edge in $F_1\/$ is $p_1\/$-parallel to $e_1\/$.

\vspace{1mm}   Likewise, $v_{2n-2}v_j \in E(F)\/$ for some $j \in \{ 2\lfloor n/2 \rfloor +1, \  2\lfloor n/2 \rfloor+3, \ \ldots, \   2n-5\}\/$. But this means that $\{v_j, v_{j+1}, \ldots, v_{2n-2}\}\/$ contains a non-crossing perfect matchings $F_2\/$ with a boundary edge $e_2\/$ and that any two edges in $F_2\/$ are $p_2\/$-parallel to $e_2\/$.

\vspace{1mm}
Since $F\/$ has only two boundary edges, we have $j = i+1\/$ (so that $i=2\lfloor n/2 \rfloor\/$ and $j = 2\lfloor n/2 \rfloor +1\/$). This proves (ii). \qed

\vspace{2mm} The following result establishes the existence of triangle-free geometric graph in regular wheel configuration arising from $k\/$  edge-disjoint  non-crossing perfect matchings each having $b\/$ boundary edges for each fixed $b \in \{2, 3\}\/$. For the case $b=3\/$, the bound on $k\/$ is best possible.

\vspace{1mm}

\begin{result} \label{b=2,3}
Let $S$ be a set of \ $2n\/$ points in regular wheel configuration in the plane.  Then there exist $k$  edge-disjoint non-crossing perfect matchings $F_1, F_2, ... , F_k$ on $S\/$ each having precisely $b\/$ boundary edges for each fixed  $b \in \{2, 3\}\/$, and $k \leq \lceil\frac{(b-1)n}{b}\rceil -(b-2)$, such that $\cup ^{k}_{i=1}{F_i} $ is a triangle-free geometric graph.  Here $n \geq b^2-1\/$.
\end{result}

\vspace{1mm}  \noindent
{\bf Proof:} (I)  $b=2$.

\vspace{1mm} Suppose $x\/$ is the center of $S\/$ and $v_0, v_1, \cdots,  v_{2n-2}\/$ are the vertices of the circle  $C\/$.  Through\ out, let $\delta\/$ take the value $1\/$  when  $n\/$ is odd, and the value $0\/$ when $n\/$ is even.    For each $i =1, 2, \ldots, k\/$, let
$$   F_i = \{ v_{i-1}x, v_{i+j-1}v_{n-j+i-\delta}, v_{n-1-\delta+i+j}v_{i-j-1} \ | \ j=1, 2, \ldots, (n-\delta)/2 \}.$$
Here the operations on the subscripts are reduced modulo $2n-1\/$. Then it is readily seen that $F_1, F_2, \ldots, F_k\/$ are $k\/$   non-crossing perfect matchings on $S\/$.

\begin{figure}[htb]
\begin{center}
\begin{tikzpicture}
%\coordinate (center) at (0,0);
%  \def\radius{2.5cm}
%   \draw (center) circle[radius=\radius];
%   \foreach \x in {0,27.6923,...,360} {
%             \filldraw[] (\x:2.5cm) circle(1pt);
%             \filldraw[] (0,0) circle(1pt);
 %              }
    %     \draw [line width=1,green](0.3,-2.5) -- (0,0);\draw [line width=1,green](0.3,2.5) -- (1.44,-2.06);\draw [line width=1,green](1.44,2.06) -- (2.22,-1.18);\draw [line width=1,green](2.22,1.18) -- (2.5,0);\draw [line width=1,green](-0.9,2.34) -- (-0.9,-2.34);\draw [line width=1,green](-2.44,0.6) -- (-2.44,-0.6);\draw [line width=1,green](-1.88,1.65) -- (-1.88,-1.65);

           %points
       \coordinate (x) at (0,0);\node at (x) {\textbullet};\node[above left] at (x) {${x}$};
       \coordinate (v0) at (0.3,-2.5);\node[below] at (v0) {${v_0}$};
       \coordinate (v6) at (0.3,2.5);\node[above] at (v6) {${v_6}$};
       \coordinate (v1) at (1.44,-2.06);\node[below] at (v1) {${v_1}$};
       \coordinate (v5) at (1.44,2.06);\node[above right] at (v5) {${v_5}$};
       \coordinate (v2) at (2.22,-1.18);\node[ right] at (v2) {${v_2}$};
       \coordinate (v4) at (2.22,1.18);\node[above right] at (v4) {${v_4}$};
       \coordinate (v3) at (2.5,0);\node[right] at (v3) {${v_3}$};
       \coordinate (v7) at (-0.9,2.34);\node[above] at (v7) {${v_7}$};
       \coordinate (v12) at (-0.9,-2.34);\node[below left] at (v12) {${v_{12}}$};
       \coordinate (v9) at (-2.44,0.6);\node[left] at (v9) {${v_{9}}$};
       \coordinate (v10) at (-2.44,-0.6);\node[left] at (v10) {${v_{10}}$};
       \coordinate (v8) at (-1.88,1.65);\node[above left] at (v8) {${v_{8}}$};
       \coordinate (v11) at (-1.88,-1.65);\node[below left] at (v11) {${v_{11}}$};

         %lines
       \draw [line width=1,green](v0) -- (x);\draw [line width=1,green](v1) -- (v6);\draw [line width=1,green](v2) -- (v5);
       \draw [line width=1,green](v3) -- (v4);\draw [line width=1,green](v7) -- (v12);\draw [line width=1,green](v8) -- (v11);
       \draw [line width=1,green](v9) -- (v10);
       \draw [line width=1,yellow](v1) -- (x);\draw [line width=1,yellow](v2) -- (v7);\draw [line width=1,yellow](v3) -- (v6);
       \draw [line width=1,yellow](v4) -- (v5);\draw [line width=1,yellow](v0) -- (v8);\draw [line width=1,yellow](v12) -- (v9);
       \draw [line width=1,yellow](v10) -- (v11);
       \draw [line width=1,blue](v2) -- (x);\draw [line width=1,blue](v3) -- (v8);\draw [line width=1,blue](v4) -- (v7);
       \draw [line width=1,blue](v5) -- (v6);\draw [line width=1,blue](v1) -- (v9);\draw [line width=1,blue](v0) -- (v10);
       \draw [line width=1,blue](v12) -- (v11);
       \draw [line width=1,red](v3) -- (x);\draw [line width=1,red](v4) -- (v9);\draw [line width=1,red](v5) -- (v8);
       \draw [line width=1,red](v6) -- (v7);\draw [line width=1,red](v2) -- (v10);\draw [line width=1,red](v1) -- (v11);
       \draw [line width=1,red](v0) -- (v12);

       \node at (x) {\textbullet}; \filldraw[black] (v0) circle(2pt);\filldraw[black] (v1) circle(2pt);\filldraw[black] (v2) circle(2pt);
       \filldraw[black] (v3) circle(2pt);
       \filldraw[black] (v4) circle(2pt);\filldraw[black] (v5) circle(2pt);\filldraw[black] (v6) circle(2pt);
       \filldraw[black] (v7) circle(2pt);\filldraw[black] (v8) circle(2pt);\filldraw[black] (v9) circle(2pt);
       \filldraw[black] (v10) circle(2pt);\filldraw[black] (v11) circle(2pt);\filldraw[black] (v12) circle(2pt);

\end{tikzpicture}
\caption{Triangle-free graph with $n=7$ vertices in regular wheel configuration.   }  \label{rwheel}
\end{center}
\end{figure}
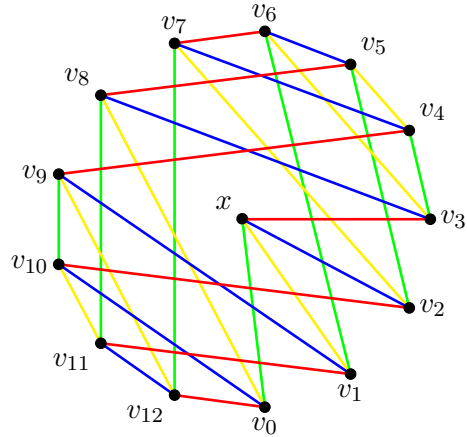

\vspace{1mm} The graph $\bigcup _{i=1} ^k F_i\/$ where $n=7\/$ is depicted in Figure \ref{rwheel}.  Note that $v_ox, v_1x, \ldots, v_{k-1}x\/$ are consecutive radial edges of $\bigcup _{i=1} ^k F_i\/$ and that $F_{i+1}\/$ is obtained from $F_i\/$ by ``rotating" the edges of $F_i\/$ with respect to the center $x\/$ of $C\/$.

\vspace{1mm} We shall show that  $G=\bigcup _{i=1} ^k F_i\/$ is triangle-free. For this purpose, note  that $k = \lceil\frac{n}{2}\rceil\/$. That is $k = (n+\delta)/2\/$.  Before proceeding, we shall take note of the following observation.

\vspace{1mm} (O1) Suppose $uv\/$ is a non-radial edge of $G\/$. Then $uv\/$ separates $V(G)-\{u, v\}\/$ into two parts $A\/$ and $A(x)\/$ with (i) $A(x)\/$ containing the center $x\/$, (ii) $|A|\/$ and $|A(x)|\/$ are both even (by Proposition \ref{s-even}), (iii) $|A| < |A(x)|\/$ and (iv)  $|A| \leq 2k - 2\delta\/$.

\vspace{1mm} Suppose $G\/$ contains a triangle $v_rv_sv_tv_r = \triangle\/$.  Then $x \not \in \{v_r, v_s, v_t\}\/$ because the neighbors $N(x)\/$  of $x\/$ are $v_0, v_1, \ldots, v_{k-1}\/$ and no two vertices in $N(x)\/$ are adjacent in $\bigcup _{i=1} ^k F_i\/$.

\vspace{2mm}
{\em Case (1):} Assume that   $x\/$ is not enclosed by $\triangle\/$.

\vspace{1mm} Suppose $\{i, j\} \subset \{r,s, t\}\/$. Let  $A_{i,j}\/$ denote the part of $V(G)-\{v_i, v_j\}\/$ separated by $v_iv_j\/$ and $A_{i,j}\/$ does not contain $x\/$. By (O1)(ii), $|A_{r,s}|, |A_{s,t}|\/$ and $|A_{t,r}|\/$ are even integers. But this is a contradiction.

\vspace{2mm}
{\em Case (2):} Assume that $x\/$ is  enclosed by $\triangle\/$.

 \vspace{1mm} Let $V(G)\/$ be partitioned into $\{x\} \cup A_0 \cup A_1 \cup A_2 \cup A_3\/$ where $A_0 = N(x)\/$ and  for each $i=1,2,3\/$,

 $A_i=\left\{
        \begin{array}{ll}
          \{v_{ik-i+1}, v_{ik-i+2}, \ldots, v_{ik+(k-i-1)} \}\/, & \hbox{when $n\/$ is odd ;} \\
          \{v_{ik}, v_{ik+1}, \ldots, v_{ik+(k-1)}\}-\{v_{4k-2} \}\/, & \hbox{otherwise.}
        \end{array}
      \right.
 $

 \vspace{1mm}
     Note that  $|A_0|=k\/$, $|A_i| = k-1\/$ when $n\/$ is odd and $|A_i| = k\/$ when $n\/$ is even.

\vspace{1mm} (O2) No two vertices in $A_j\/$ are adjacent for any $j \in \{1, 3\}\/$ (by the definition of $G\/$).

\vspace{1mm} (O3) Suppose $uv\/$ is a non-radial edge that separates $V(G) -\{u, v\}\/$ into $A\/$ and $A(x)\/$. Then $uv \in E(F_i)\/$ for some non-crossing perfect matching $F_i\/$ in $G\/$.  As such,  $|A \cap A_j| \leq |A|/2\/$ for any  $j \in \{1, 3\}\/$.

 \vspace{1mm} To see that (O3) is true, assume the contrary that $|A \cap A_j|  > |A|/2\/$. Then  either there is some boundary edge in $A \cap A_j\/$ belonging to $F_i\/$, or else some edge in $F_i\/$ incident to a vertex in $A \cap A_j\/$ is crossing with some other edges of $F_i\/$. Either case is a contradiction.

\vspace{1mm}  Suppose the triangle $\triangle\/$ is such that $v_r \in A_1\/$, $v_s \in A_2\/$ and $v_t \in A_4\/$. Then the edge $v_sv_t\/$ separates $V(G) - \{v_s, v_t\}\/$ into $A\/$ and $A_x\/$ with $A_3 \subset A\/$. But this is impossible by (O1)(iv) and (O3).

\vspace{1mm} Hence, by symmetry, we may assume without loss of generality that $v_r \in A_1\/$, $v_s \in A_2\/$ and $v_t \in A_3\/$. Then with the same notations as in Case (1), we have $|A_{r, t}| \leq 2k -2\delta\/$ (by (O1)(iv)) and $|A_{r,s} \cup A_{s,t}| \geq 2k-1\/$ (since $|C| = 4k -2\delta-1\/$). But this implies that either $|A_1 \cap A_{r,s}| > |A_{r,s}|/2\/$ or else $|A_3 \cap A_{s,t}| > |A_{s,t}|/2\/$. Either case contradicts (O3).

\vspace{1mm} This completes the proof for Case (2).

\vspace{2mm} (II)  $b=3$.

\vspace{1mm} Suppose $x$ is the center of $S$ and $v_0, v_1, \cdots,  v_{2n-2}\/$ are the vertices of the circle  $C\/$.  We consider three cases $n =3m+\ell \geq 8\/$ where $\ell \in \{0, 1, 2\}\/$.

\vspace{1mm} When  $n=3m$,  we let

\vspace{1mm}
$F_1=\{v_0x,  v_{3m-i} v_{3m-1+i},  v_jv_{2m+1-j},  v_{-j} v_{j-2m-1} \ | \ i=1,2, ...,m-1, \  j=1,2,...,m \}$.

\vspace{1mm}
 When $n=3m+1$,  we  let

\vspace{1mm}
$F_1=\{v_ox,  v_{3m-1-i} v_{3m-2+i},  v_j v_{2m-1-j},  v_{-r} v_{r-2m-3}  \ | \ i=1,2, ...,m, \ j=1,2,...,m-1, \  r=1,2,...,m+1 \}$.

\vspace{1mm}
 When $n=3m+2$,  we let

 \vspace{1mm}
$F_1=\{v_0x,  v_{3m+2-i} v_{3m+1+i},  v_j v_{2m+1-j},  v_{-j} v_{j-2m-1}  \ | \  i=1,2, \dots,m+1, \ j=1,2,...,m \}$. \

\vspace{2mm}
Now for each $p=2,...,k$, let  $F_p$ be obtained from $F_1$ by replacing each $v_s$ in $F_1$ by $v_{s+2p-2}$.\

\begin{figure}
\resizebox{15cm}{!}{
\begin{minipage}{.45\textwidth}
\begin{tikzpicture}
\coordinate (center) at (0,0);
  \def\radius{2.5cm}
 %  \draw (center) circle[radius=\radius];
   \foreach \x in {0,21.1764,...,360} {
             \filldraw[] (\x:2.5cm) circle(1pt);
             \filldraw[] (0,0) circle(1pt);
               }

       \coordinate (x) at (0,0);\node at (x) {\textbullet};\node[above left] at (x) {${x}$};
       \coordinate (v0) at (-0.7,-2.4);\node at (v0) {\textbullet};\node[below] at (v0) {${v_0}$};
       \coordinate (v1) at (0.25,-2.48);\node at (v1) {\textbullet};\node[below] at (v1) {${v_1}$};
       \coordinate (v2) at (1.1,-2.2);\node at (v2) {\textbullet};\node[below right] at (v2) {${v_2}$};
       \coordinate (v3) at (1.86,-1.7);\node at (v3) {\textbullet};\node[below right] at (v3) {${v_3}$};
       \coordinate (v4) at (2.36,-.9);\node at (v4) {\textbullet};\node[below right] at (v4) {${v_4}$};
       \coordinate (v5) at (2.5,0);\node at (v5) {\textbullet};\node[below right] at (v5) {${v_5}$};
       \coordinate (v6) at (2.36,0.9);\node at (v6) {\textbullet};\node[right] at (v6) {${v_6}$};
       \coordinate (v7) at (1.86,1.7);\node at (v7) {\textbullet};\node[above right] at (v7) {${v_7}$};
       \coordinate (v8) at (1.1,2.2);\node at (v8) {\textbullet};\node[above right] at (v8) {${v_8}$};
       \coordinate (v9) at (0.25,2.48);\node at (v9) {\textbullet};\node[above] at (v9) {${v_9}$};
       \coordinate (v10) at (-0.7,2.4);\node at (v10) {\textbullet};\node[above left] at (v10) {${v_{10}}$};
       \coordinate (v11) at (-1.5,2);\node at (v11) {\textbullet};\node[above left] at (v11) {${v_{11}}$};
       \coordinate (v12) at (-2.1,1.3);\node at (v12) {\textbullet};\node[above left] at (v12) {${v_{12}}$};
       \coordinate (v13) at (-2.46,0.5);\node at (v13) {\textbullet};\node[above left] at (v13) {${v_{13}}$};
       \coordinate (v14) at (-2.46,-0.5);\node at (v14) {\textbullet};\node[below left] at (v14) {${v_{14}}$};
       \coordinate (v15) at (-2.1,-1.3);\node at (v15) {\textbullet};\node[below left] at (v15) {${v_{15}}$};
       \coordinate (v16) at (-1.5,-2);\node at (v16) {\textbullet};\node[below left] at (v16) {${v_{16}}$};

         %lines
       \draw [line width=1,green](v0) -- (x);\draw [line width=1,green](v1) -- (v6);\draw [line width=1,green](v2) -- (v5);
       \draw [line width=1,green](v3) -- (v4);\draw [line width=1,green](v16) -- (v11);\draw [line width=1,green](v15) -- (v12);
       \draw [line width=1,green](v14) -- (v13);\draw [line width=1,green](v7) -- (v10);\draw [line width=1,green](v8) -- (v9);
       \draw [line width=1,blue](v2) -- (x);\draw [line width=1,blue](v3) -- (v8);\draw [line width=1,blue](v4) -- (v7);
       \draw [line width=1,blue](v5) -- (v6);\draw [line width=1,blue](v1) -- (v13);\draw [line width=1,blue](v0) -- (v14);
       \draw [line width=1,blue](v16) -- (v15);\draw [line width=1,blue](v9) -- (v12);\draw [line width=1,blue](v10) -- (v11);
       \draw [line width=1,red](v4) -- (x);\draw [line width=1,red](v5) -- (v10);\draw [line width=1,red](v6) -- (v9);
       \draw [line width=1,red](v7) -- (v8);\draw [line width=1,red](v3) -- (v15);\draw [line width=1,red](v2) -- (v16);
       \draw [line width=1,red](v1) -- (v0);\draw [line width=1,red](v11) -- (v14);\draw [line width=1,red](v12) -- (v13);
       \draw [line width=1,yellow](v6) -- (x);\draw [line width=1,yellow](v7) -- (v12);\draw [line width=1,yellow](v8) -- (v11);
       \draw [line width=1,yellow](v9) -- (v10);\draw [line width=1,yellow](v5) -- (v0);\draw [line width=1,yellow](v4) -- (v1);
       \draw [line width=1,yellow](v3) -- (v2);\draw [line width=1,yellow](v13) -- (v16);\draw [line width=1,yellow](v14) -- (v15);
       \draw [line width=1,brown](v8) -- (x);\draw [line width=1,brown](v9) -- (v14);\draw [line width=1,brown](v10) -- (v13);
       \draw [line width=1,brown](v11) -- (v12);\draw [line width=1,brown](v7) -- (v2);\draw [line width=1,brown](v6) -- (v3);
       \draw [line width=1,brown](v5) -- (v4);\draw [line width=1,brown](v15) -- (v1);\draw [line width=1,brown](v16) -- (v0);

       \filldraw[black] (x) circle(2pt);\filldraw[black] (v0) circle(2pt);\filldraw[black] (v1) circle(2pt);
       \filldraw[black] (v2) circle(2pt);\filldraw[black] (v3) circle(2pt);\filldraw[black] (v4) circle(2pt);
       \filldraw[black] (v5) circle(2pt);\filldraw[black] (v6) circle(2pt);\filldraw[black] (v7) circle(2pt);
       \filldraw[black] (v8) circle(2pt);\filldraw[black] (v9) circle(2pt);\filldraw[black] (v10) circle(2pt);
       \filldraw[black] (v11) circle(2pt);\filldraw[black] (v12) circle(2pt);\filldraw[black] (v13) circle(2pt);
       \filldraw[black] (v14) circle(2pt);\filldraw[black] (v15) circle(2pt);\filldraw[black] (v16) circle(2pt);

\end{tikzpicture}
\centering

(a) $n=3m=9$
\end{minipage}
\begin{minipage}{.45\textwidth}
\begin{tikzpicture}
%\coordinate (center) at (0,0);
%  \def\radius{2.5cm}
 %  \draw (center) circle[radius=\radius];
 %  \foreach \x in {0,24,...,360} {
  %           \filldraw[] (\x:2.5cm) circle(1pt);
   %          \filldraw[] (0,0) circle(1pt);
    %           }
    %     \draw [line width=1,green](-0.26,-2.5) -- (0,0);\draw [line width=1,green](-0.26,2.5) -- (0.78,-2.4);\draw [line width=1,green](0.78,2.4) -- (1.68,-1.9);\draw [line width=1,green](1.68,1.9) -- (2.28,-1);\draw [line width=1,green](-1.25,2.2) -- (-1.25,-2.2);\draw [line width=1,green](-2.02,1.46) -- (-2.02,-1.46);\draw [line width=1,green](-2.44,0.54) -- (-2.44,-0.54);
           %points
       \coordinate (x) at (0,0);\node at (x) {\textbullet};\node[above left] at (x) {${x}$};
       \coordinate (v0) at (-0.26,-2.5);\filldraw[black] (v0) circle(2pt);\node[below] at (v0) {${v_0}$};
       \coordinate (v8) at (-0.26,2.5);\filldraw[black] (v8) circle(2pt);\node[above] at (v8) {${v_8}$};
       \coordinate (v1) at (0.78,-2.4);\filldraw[black] (v1) circle(2pt);\node[below right] at (v1) {${v_1}$};
       \coordinate (v7) at (0.78,2.4);\filldraw[black] (v7) circle(2pt);\node[above] at (v7) {${v_7}$};
       \coordinate (v2) at (1.68,-1.9);\filldraw[black] (v2) circle(2pt);\node[below right] at (v2) {${v_2}$};
       \coordinate (v6) at (1.68,1.9);\filldraw[black] (v6) circle(2pt);\node[above right] at (v6) {${v_6}$};
       \coordinate (v3) at (2.28,-1);\filldraw[black] (v3) circle(2pt);\node[below right] at (v3) {${v_3}$};
       \coordinate (v5) at (2.28,1);\filldraw[black] (v5) circle(2pt);\node[above right] at (v5) {${v_5}$};
       \coordinate (v4) at (2.5,0);\filldraw[black] (v4) circle(2pt);\node[right] at (v4) {${v_4}$};
       \coordinate (v9) at (-1.25,2.2);\filldraw[black] (v9) circle(2pt);\node[above left] at (v9) {${v_9}$};
       \coordinate (v14) at (-1.25,-2.2);\filldraw[black] (v14) circle(2pt);\node[below left] at (v14) {${v_{14}}$};
       \coordinate (v10) at (-2.02,1.46);\filldraw[black] (v10) circle(2pt);\node[above left] at (v10) {${v_{10}}$};
       \coordinate (v13) at (-2.02,-1.46);\filldraw[black] (v13) circle(2pt);\node[below left] at (v13) {${v_{13}}$};
      \coordinate (v11) at (-2.44,0.54);\filldraw[black] (v11) circle(2pt);\node[left] at (v11) {${v_{11}}$};
     \coordinate (v12) at (-2.44,-0.54);\filldraw[black] (v12) circle(2pt);\node[left] at (v12) {${v_{12}}$};

         %lines
              \draw [line width=1,green](v0) -- (x);\draw [line width=1,green](v1) -- (v4);\draw [line width=1,green](v2) -- (v3);
       \draw [line width=1,green](v14) -- (v11);\draw [line width=1,green](v13) -- (v12);\draw [line width=1,green](v5) -- (v10);
       \draw [line width=1,green](v6) -- (v9);\draw [line width=1,green](v7) -- (v8);
       \draw [line width=1,blue](v2) -- (x);\draw [line width=1,blue](v3) -- (v6);\draw [line width=1,blue](v4) -- (v5);
       \draw [line width=1,blue](v1) -- (v13);\draw [line width=1,blue](v0) -- (v14);\draw [line width=1,blue](v7) -- (v12);
       \draw [line width=1,blue](v8) -- (v11);\draw [line width=1,blue](v9) -- (v10);
       \draw [line width=1,red](v4) -- (x);\draw [line width=1,red](v5) -- (v8);\draw [line width=1,red](v6) -- (v7);
       \draw [line width=1,red](v3) -- (v0);\draw [line width=1,red](v2) -- (v1);\draw [line width=1,red](v9) -- (v14);
       \draw [line width=1,red](v10) -- (v13);\draw [line width=1,red](v11) -- (v12);
       \draw [line width=1,yellow](v6) -- (x);\draw [line width=1,yellow](v7) -- (v10);\draw [line width=1,yellow](v8) -- (v9);
       \draw [line width=1,yellow](v5) -- (v2);\draw [line width=1,yellow](v4) -- (v3);\draw [line width=1,yellow](v11) -- (v1);
       \draw [line width=1,yellow](v12) -- (v0);\draw [line width=1,yellow](v13) -- (v14);
       \draw [line width=1,brown](v8) -- (x);\draw [line width=1,brown](v9) -- (v12);\draw [line width=1,brown](v10) -- (v11);
       \draw [line width=1,brown](v7) -- (v4);\draw [line width=1,brown](v6) -- (v5);\draw [line width=1,brown](v13) -- (v3);
       \draw [line width=1,brown](v14) -- (v2);\draw [line width=1,brown](v0) -- (v1);

       \filldraw[black] (x) circle(2pt);\filldraw[black] (v0) circle(2pt);\filldraw[black] (v1) circle(2pt);
       \filldraw[black] (v2) circle(2pt);\filldraw[black] (v3) circle(2pt);\filldraw[black] (v4) circle(2pt);
       \filldraw[black] (v5) circle(2pt);\filldraw[black] (v6) circle(2pt);\filldraw[black] (v7) circle(2pt);
       \filldraw[black] (v8) circle(2pt);\filldraw[black] (v9) circle(2pt);\filldraw[black] (v10) circle(2pt);
       \filldraw[black] (v11) circle(2pt);\filldraw[black] (v12) circle(2pt);\filldraw[black] (v13) circle(2pt);
       \filldraw[black] (v14) circle(2pt);

\end{tikzpicture}
\centering

(b) $n=3m+2=8$
\end{minipage}
\begin{minipage}{.45\textwidth}
\begin{tikzpicture}
\coordinate (center) at (0,0);
  \def\radius{2.5cm}
%   \draw (center) circle[radius=\radius];
   \foreach \x in {0,18.947,...,360} {
             \filldraw[] (\x:2.5cm) circle(1pt);
             \filldraw[] (0,0) circle(1pt);
                }
            %points
       \coordinate (x) at (0,0);\node at (x) {\textbullet};\node[above left] at (x) {${x}$};
       \coordinate (v0) at (-0.2,-2.5);\node at (v0) {\textbullet};\node[below] at (v0) {${v_0}$};
       \coordinate (v1) at (0.6,-2.4);\node at (v1) {\textbullet};\node[below right] at (v1) {${v_1}$};
       \coordinate (v2) at (1.4,-2.1);\node at (v2) {\textbullet};\node[below right] at (v2) {${v_2}$};
       \coordinate (v3) at (2,-1.54);\node at (v3) {\textbullet};\node[below right] at (v3) {${v_3}$};
       \coordinate (v4) at (2.38,-0.8);\node at (v4) {\textbullet};\node[below right] at (v4) {${v_4}$};
       \coordinate (v5) at (2.5,0);\node at (v5) {\textbullet};\node[right] at (v5) {${v_5}$};
       \coordinate (v6) at (2.38,0.8);\node at (v6) {\textbullet};\node[above right] at (v6) {${v_6}$};
       \coordinate (v7) at (2,1.54);\node at (v7) {\textbullet};\node[above right] at (v7) {${v_7}$};
       \coordinate (v8) at (1.4,2.1);\node at (v8) {\textbullet};\node[above right] at (v8) {${v_8}$};
       \coordinate (v9) at (0.6,2.4);\node at (v9) {\textbullet};\node[above right] at (v9) {${v_9}$};
       \coordinate (v10) at (-0.2,2.5);\node at (v10) {\textbullet};\node[above] at (v10) {${v_{10}}$};
       \coordinate (v11) at (-1,2.34);\node at (v11) {\textbullet};\node[above left] at (v11) {${v_{11}}$};
       \coordinate (v12) at (-1.72,1.86);\node at (v12) {\textbullet};\node[above left] at (v12) {${v_{12}}$};
       \coordinate (v13) at (-2.2,1.2);\node at (v13) {\textbullet};\node[above left] at (v13) {${v_{13}}$};
       \coordinate (v14) at (-2.46,.4);\node at (v14) {\textbullet};\node[left] at (v14) {${v_{14}}$};
       \coordinate (v15) at (-2.46,-.4);\node at (v15) {\textbullet};\node[left] at (v15) {${v_{15}}$};
       \coordinate (v16) at (-2.2,-1.2);\node at (v16) {\textbullet};\node[left] at (v16) {${v_{16}}$};
       \coordinate (v17) at (-1.72,-1.86);\node at (v17) {\textbullet};\node[below left] at (v17) {${v_{17}}$};
       \coordinate (v18) at (-1,-2.34);\node at (v18) {\textbullet};\node[below left] at (v18) {${v_{18}}$};

         %lines
        \draw [line width=1,green](v0) -- (x);\draw [line width=1,green](v1) -- (v4);\draw [line width=1,green](v2) -- (v3);
       \draw [line width=1,green](v5) -- (v10);\draw [line width=1,green](v6) -- (v9);\draw [line width=1,green](v7) -- (v8);
       \draw [line width=1,green](v18) -- (v11);\draw [line width=1,green](v17) -- (v12);\draw [line width=1,green](v16) -- (v13);
       \draw [line width=1,green](v15) -- (v14);
       \draw [line width=1,blue](v2) -- (x);\draw [line width=1,blue](v3) -- (v6);\draw [line width=1,blue](v4) -- (v5);
       \draw [line width=1,blue](v7) -- (v12);\draw [line width=1,blue](v8) -- (v11);\draw [line width=1,blue](v9) -- (v10);
       \draw [line width=1,blue](v1) -- (v13);\draw [line width=1,blue](v0) -- (v14);\draw [line width=1,blue](v18) -- (v15);
       \draw [line width=1,blue](v17) -- (v16);
       \draw [line width=1,red](v4) -- (x);\draw [line width=1,red](v5) -- (v8);\draw [line width=1,red](v6) -- (v7);
       \draw [line width=1,red](v9) -- (v14);\draw [line width=1,red](v10) -- (v13);\draw [line width=1,red](v11) -- (v12);
       \draw [line width=1,red](v3) -- (v15);\draw [line width=1,red](v2) -- (v16);\draw [line width=1,red](v1) -- (v17);
       \draw [line width=1,red](v0) -- (v18);
       \draw [line width=1,brown](v6) -- (x);\draw [line width=1,brown](v7) -- (v10);\draw [line width=1,brown](v8) -- (v9);
       \draw [line width=1,brown](v11) -- (v16);\draw [line width=1,brown](v12) -- (v15);\draw [line width=1,brown](v13) -- (v14);
       \draw [line width=1,brown](v5) -- (v17);\draw [line width=1,brown](v4) -- (v18);\draw [line width=1,brown](v3) -- (v0);
       \draw [line width=1,brown](v2) -- (v1);
       \draw [line width=1,black](v8) -- (x);\draw [line width=1,black](v9) -- (v12);\draw [line width=1,black](v10) -- (v11);
       \draw [line width=1,black](v13) -- (v18);\draw [line width=1,black](v14) -- (v17);\draw [line width=1,black](v15) -- (v16);
       \draw [line width=1,black](v7) -- (v0);\draw [line width=1,black](v6) -- (v1);\draw [line width=1,black](v5) -- (v2);
       \draw [line width=1,black](v4) -- (v3);
       \draw [line width=1,yellow](v10) -- (x);\draw [line width=1,yellow](v11) -- (v14);\draw [line width=1,yellow](v12) -- (v13);
       \draw [line width=1,yellow](v15) -- (v1);\draw [line width=1,yellow](v16) -- (v0);\draw [line width=1,yellow](v17) -- (v18);
       \draw [line width=1,yellow](v9) -- (v2);\draw [line width=1,yellow](v8) -- (v3);\draw [line width=1,yellow](v7) -- (v4);
       \draw [line width=1,yellow](v6) -- (v5);

       \filldraw[black] (x) circle(2pt);\filldraw[black] (v0) circle(2pt);\filldraw[black] (v1) circle(2pt);
       \filldraw[black] (v2) circle(2pt);\filldraw[black] (v3) circle(2pt);\filldraw[black] (v4) circle(2pt);
       \filldraw[black] (v5) circle(2pt);\filldraw[black] (v6) circle(2pt);\filldraw[black] (v7) circle(2pt);
       \filldraw[black] (v8) circle(2pt);\filldraw[black] (v9) circle(2pt);\filldraw[black] (v10) circle(2pt);
       \filldraw[black] (v11) circle(2pt);\filldraw[black] (v12) circle(2pt);\filldraw[black] (v13) circle(2pt);
       \filldraw[black] (v14) circle(2pt);\filldraw[black] (v15) circle(2pt);\filldraw[black] (v16) circle(2pt);
       \filldraw[black] (v17) circle(2pt);\filldraw[black] (v18) circle(2pt);

\end{tikzpicture}
\centering

(c) $n=3m+1=10$
\end{minipage}
}
\caption{ Triangle-free geometric graphs with $b=3\/$ boundary edges.    \label{2wheel}  }
\end{figure}
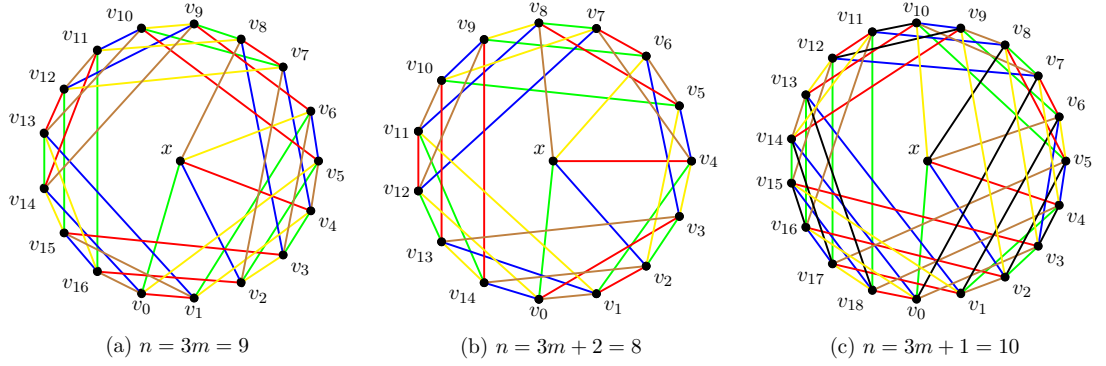

%example
\vspace{1mm}
The above constructions are illustrated in Figure \ref{2wheel}(a), (b) and (c) for the cases $n=9\/$, $n=8\/$ and $10\/$ respectively.

\vspace{1mm} Now we shall show that  $G=\bigcup _{i=1} ^k F_i\/$ is triangle-free. For this purpose, we may assume that $k = \lceil \frac{2n}{3} \rceil -1\/$.

\vspace{1mm}  First we observe that, similar to part (I), if $G\/$ contains a triangle $v_rv_sv_tv_r = \triangle\/$ where $r<s<t$, then $x \not \in \{v_r, v_s, v_t\}\/$ because the neighbors $N(x)\/$  of $x\/$ are $v_0, v_2, \ldots, v_{2k-2}\/$ and no two vertices in $N(x)\/$ are adjacent in $\bigcup _{i=1} ^k F_i\/$.

\vspace{1mm} Next we take note of the following observation which is similar as in part (I).

\vspace{1mm} (O4) Suppose $v_iv_j\/$ is a non-radial edge of $G\/$. Then $v_iv_j\/$ separates $V(G)-\{v_i, v_j\}\/$ into two parts $A_{i,j}\/$ and $A_{i,j}(x)\/$ with (i) $A_{i,j}(x)\/$ containing the center $x\/$, (ii) $|A_{i,j}|\/$ and $|A_{i,j}(x)|\/$ are both even (by Proposition \ref{s-even}), (iii) $|A_{i,j}| < |A_{i,j}(x)|\/$ and (iv)  $|A_{i,j}| \leq 2m-2+\ell\/$ where $\ell=0,2$ and $|A_{i,j}| \leq 2m\/$ where $\ell=1$ (recall that $n =3m + \ell$).

\vspace{2mm}
{\em Case (1):} Assume that   $x\/$ is not enclosed by $\triangle\/$.

\vspace{1mm} By relabeling $r, s, t\/$, if necessary, we may assume without loss of generality that $v_r, v_s\/$ separates  $G -\{v_r, v_s\}\/$ into $A_{r,s}\/$ and $A_{r,s}(x)$ so that $v_t \in A_{r,s}\/$. Then by (O1)(ii), $|A_{r,s}(x)|, |A_{s,t}|\/$ and $|A_{t,r}|\/$ are even integers. But this is a contradiction.

\vspace{2mm}
{\em Case (2):} Assume that $x\/$ is  enclosed by $\triangle\/$.

\vspace{1mm}
 (i) Suppose  only one of the  $A_{i,j}$'s attains the maximum number of vertices.

\vspace{1mm} When $\ell \in \{ 0, 2\}$, we have  $|A_{r,s}|+|A_{s,t}|+|A_{t,r}| \leq 6m -10 +3\ell <6m+2\ell-4=2n-4$ which is a contradiction. Also when $\ell=1$, we have  $|A_{r,s}|+|A_{s,t}|+|A_{t,r}|<6m-2=2n-4$ and is  again a contradiction.

\vspace{2mm}

(ii) Suppose only two parts of the $A_{i,j}$'s attain the maximum number of vertices.

\vspace{1mm} When $\ell \in \{ 0, 2\}\/$, we have  $|A_{r,s}|+|A_{s,t}|+|A_{t,r}| \leq 6m+3\ell-8  <6m+2\ell-4=2n-4$ which is a contradiction.

\vspace{1mm} Hence assume that $\ell=1$.

\vspace{1mm} Suppose  $A_{r,s}$ and $A_{t,r}\/$ attain the  maximum number of vertices. Assume further $v_rv_s\in F'$ and $v_rv_t\in F''$, where $F'$ and $F''$ are non-crossing perfect matchings.

\vspace{1mm} It is clear that $v_{s+1}x \in F'$ and $v_{r+1}x\in F''$ (since $v_{r-1}x\notin F'$ and $v_{t-1}x\notin F''$ by definition of $G\/$) and this  implies that the number of vertices in $v_{r+2},v_{r+3},...,v_{s}$ is even. This  contradicts the definition of $G\/$ (since the indices of the radial vertices are of the same parity).

\vspace{1mm}
(iii) Suppose all the  $A_{i,j}$'s  attain the maximum number of vertices.

\vspace{1mm} Note that if $v_iv_j\in F'\/$ for some non-crossing perfect matching $F'\/$ of $G\/$ where $i<j$.  Then either $v_{i-1}x$ or $v_{j+1}x$ is a radial edge of $F'\/$ (but not both by definition of $G\/$). % By symmetry, there are only two choices to consider.

 \vspace{1mm} Hence each $A_{i,j}\/$ contains  only one radial vertex.     Without loss of generality assume that $v_{r+1}x$, $v_{s+1}x$ and $v_{t+1}x$ are radial edges  with $v_{r+1} \in A_{r,s}$, $v_{s+1} \in A_{s,t}$ and $v_{t+1} \in A_{t,r}$.  This implies that the number of vertices in $v_{i+1}, v_{i+2},..., v_{j-1} \in A_{i,j}$ is even (by Proposition \ref{s-even}) for any $i, j \in \{r, s, t\}\/$. This contradicts the definition of $G\/$ (since the indices of the radial vertices are of the same parity).

%\vspace{1mm} (b) Some  $A_{i,j}\/$ contains  two radial vertices.

%\vspace{1mm} Without loss of generality assume that $v_{r+1}x$ and $v_{s-1}x$ are radial edges with $v_{r+1}, v_{s-1} \in A_{r,s}$.   But  implies that the number of vertices in $v_{r+2},v_{r+2},...,v_{s-2}$ is even which contradicts the definition of $G\/$ (since the indices of the radial vertices are of the same parity).

\vspace{1mm} This completes the proof.  \qed

\vspace{1mm}
The next result shows that, when $n \geq 3\/$ is odd and  if every non-crossing perfect matching has only two boundary edges, the geometric graph obtained with $\lceil n/2 \rceil\/$ perfect matchings (that satisfies the conditions of Theorem \ref{b=2,3}) is unique. Note that, in this case,  the  set of radial edges are in consecutive order (i.e, in the form  $xv_0, xv_1, \ldots , xv_{k-1}\/$).

\begin{result}\label{w-oddn2}
 Let $S\/$ be a set of $2n\/$ points in regular wheel configuration  in the plane where $n \geq 3\/$ is odd.
 Suppose $F_1, F_2, \ldots, F_k\/$ are  $k\/$  edge-disjoint non-crossing perfect matchings on $S\/$ such that $G=\bigcup _{i=1} ^k F_i\/$ is a triangle-free geometric graph with  $k = \lceil n/2 \rceil\/$. Then each  non-crossing perfect matching $F_i\/$ has exactly two boundary edges if and only if all radial edges in $G\/$ are in consecutive order.
\end{result}

\vspace{1mm}  \noindent
{\bf Proof:}
 We first suppose that each $F_i\/$ has only two boundary edges. Assume without loss of generality that $F_1 = \{ v_0x, v_jv_{n-j}, v_{n-1+j}v_{-j} \ | \ j=1, 2, \ldots, (n-1)/2 \} .\/$

%\vspace{2mm} {\bf ???} Note that if the radial edge of $F_i\/$ is $v_px\/$, then the boundary edges of $F_i\/$ are $v_{k-1+p}v_{k+p}\/$ and  $v_{3k-3+p}v_{3k-2+p}\/$.

\vspace{2mm} Let $F\/$ be a non-crossing perfect matching in $G\/$ and $v_px\/$ is the radial edge in $F\/$. Suppose  $F'\/$ is a non-crossing perfect matching in $G\/$  different from $F\/$ and $v_qx\/$ is the radial edge of $F'\/$, then $q\/$ is neither $p+2k-2\/$ nor $p+2k-1\/$; otherwise $E(F') \cap E(F) \neq \emptyset\/$. That is

\vspace{1mm} (i) if $v_px \/$ is a radial edge of $G\/$, then neither  $v_{p+2k-2}x\/$ nor $v_{p+2k-1}x\/$ can be a radial edge of $G\/$.

\vspace{2mm}
Assume on the contrary that not all radial edges of $G\/$ are in consecutive order. By symmetry, we may assume without loss of generality that there exists $q \in \{k, k+1, \ldots, 2k-3\}\/$ such that $v_qx\/$ is the radial edge of a non-crossing perfect matching $F'$ in $G\/$.

\vspace{2mm}
 We claim that, for each $q \in \{k, k+1, \ldots, 2k-3\}\/$ where $v_qx\/$ and $v_0x\/$ are both radial edges of $G\/$, there exists at most $k-3\/$ other vertices   which can be used as radial vertices of $G\/$, thereby establishing a contradiction (because $G\/$ has precisely $k\/$ radial edges).\

\vspace{2mm} First, we show that the claim is true if $q=k\/$.

\vspace{1mm}  Since $v_kx\/$ is the radial edge of $F'\/$, we have
$F' = \{ v_kx, v_{k+j}v_{n-j+k}, v_{n-1+k+j}v_{k-j} \ | \ j=1, 2, \ldots, (n-1)/2 \}.$ \ Further,  by (i), we have

\vspace{1mm} (r1) $v_j\/$ is not  a radial vertex  of $G\/$ for any $j \in \{2k-2, 2k-1, 3k-2, 3k-1\}\/$.

\vspace{2mm} Observe that,  for any $l= 2, 3, \ldots, k-2\/$, $v_{k-l}x\/$ is not a radial edge of $G\/$ otherwise $v_{3k-l-1}v_{k-l-1} \in E(G)\/$ which yields a triangle $v_{k+l}v_{3k-l-1}v_{k-l-1}  v_{k+l}\/$ in $G\/$, a contradiction.

\vspace{2mm}
Clearly, for $i \in \{1, k-1\}\/$,  $v_ix\/$ is not a radial edge of $G\/$ (otherwise $xv_{i-1}v_ix\/$ will be a triangle in $G\/$). That means

\vspace{1mm} (r2) $v_{k-l}\/$ is not a radial vertex of $G\/$ for every $l = 1, 2, \ldots, k-1\/$.

\vspace{2mm}
For any $l=0, 1, \ldots, k-3\/$, $v_{2k+l}x\/$ is not a radial edge of $G\/$ otherwise $v_{2k+l-1}v_{l+2} \in E(G)\/$ which yields a triangle $v_{2k+l-1}v_{l+2} v_{-l-1}v_{2k+l-1}\/$ in $G\/$, a contradiction. That means

\vspace{1mm} (r3)  $v_{2k+l}\/$ is not a radial vertex of $G\/$  for any $l = 0, 1, \ldots, k-3\/$.

\vspace{2mm}
By (i), we have
%Further, for any $m=1, 2, \ldots, k-3\/$, if $v_{k+m}x\/$ and $v_{3k+m-1}x\/$ are radial edges of $F_3$ and $F_4$ in $G\/$. Then $E(F_3) \cap E(F_4) \neq \emptyset\/$ (since $F_i\/$ has only two boundary edges) . That means\\

\vspace{1mm} (r4)  for each  $l =1, 2, \ldots, k-3$,  $v_{k+l}\/$ and $v_{3k+l-1}\/$  cannot be both radial vertices of $G\/$.

\vspace{2mm} It  follows from  (r1), (r2), (r3) and (r4) that  the claim is true for $q=k\/$.

\vspace{2mm} Now let $q = k+r\/$ for some $r \in \{1, 2,\ldots,k-3\}\/$ and assume that $v_qx\/$ is a radial edge of $G\/$. Also assume (it has been shown) that $v_{k+i}x\/$ is not a radial edge of $G\/$ for any $i=0, 1, \ldots, r-1\/$.
With this assumption, it  implies that  $v_ix$ is not  a radial edge of $G\/$ for any $1 \leq i \leq r\/$ otherwise having both $v_ix\/$ and  $v_{k+r}x\/$ as radial edges of $G\/$ would mean that  $v_0x\/$ and $v_{k+r-i}x\/$ are both radial edges of $G\/$, a contradiction (to the induction hypothesis).

%Suppose $v_{q+1}x=v_{k+r}\/$ is a radial edge of $G\/$.
\vspace{2mm} As such, we have

\vspace{1mm} $F'=\{v_{k+r}x,v_{k+r+j}v_{3k-1+r-j},v_{3k-2+r+j}v_{k+r-j}\ | \ j=1,2,...,(n-1)/2 \}\/$.

Further by (i), we have\

\vspace{1mm}
(w1) for any $s\in \{2k-2,2k-1,3k+r-2,3k+r-1\}$, $v_s$ is not a radial vertex of $G\/$.

\vspace{1mm} Also by (i),

\vspace{1mm} (w2) for each $s=r+1,r+2,...,k-3$, $v_{k+s}$ and $v_{3k+s-1}$ cannot be both radial vertices of $G\/$.

\vspace{1mm}
 Now for any $1\leq s\leq \lfloor(k-r-1)/2\rfloor$, $v_{r+s}x$ is not a radial edge of $G\/$ otherwise $v_{2r+s+1}v_{2k+s-2} \in E(G)\/$ which yields a triangle $v_{-s}v_{2r+s+1}v_{2k+s-2} v_{-s}\/$ in $G\/$, a contradiction.

\vspace{1mm}
Also, for any $\lfloor(k-r-1)/2\rfloor+1\leq s\leq k-r-1$, $v_{r+s}x$ is not a radial edge of $G\/$ otherwise $v_{2k+2r+s-1}v_{s-1} \in E(G)\/$ which yields a triangle $v_{s-1}v_{2k-s}v_{2k+2r+s-1}v_{s-1}\/$ in $G\/$, a contradiction. That is

\vspace{1mm} (w3) for any $j=r+1,r+2,...,k-1$,   $v_j\/$ is not a radial vertex of $G\/$.

\vspace{1mm}
For any $i=r,r+1,...,k-3$, $v_{2k+i}x\/$ is  not a radial edge of $G\/$, otherwise $v_{2k+i-r-1}v_{r+i+2} \linebreak \in E(G)\/$ which yields a triangle $v_{r-i-1}v_{2k+i-r-1}v_{r+i+2}v_{r-i-1}\/$ in $G\/$, a contradiction. That is

%For any $h=0,1, \ldots,k-r-3$. $v_{2k+r+h}x$ is not a radial edge of $G\/$, otherwise $v_{2k+h-1}v_{h+2r+2} \in E(G)\/$ which yields a triangle $v_{-h-1}v_{2k+h-1}v_{h+2r+2}v_{-h-1}\/$ in $G\/$, a contradiction. That means\

\vspace{1mm} (w4) for any $i=r,r+1,...,k-3$,  $v_{2k+i}\/$ is not a radial vertex of $G\/$.

\vspace{2mm}Let $j\in \{1,2, \ldots,r \}$.
 Suppose $v_{2k+j-1}x$ and $v_{3k+ j-3}x$ are radial edges of $G$. Then $v_{2j}v_{2k-1}\in E(G)\/$ which yields a triangle $v_{-1} v_{2j}v_{2k-1}v_{-1}\/$ in $G\/$, a contradiction. That is,

 \vspace{1mm} (w5)  $v_{2k+j-1}\/$ and $v_{3k+j-3}$ cannot be both radial vertices of $G\/$ for any $j=1, \ldots,r $.\

\vspace{2mm} By (w1),(w2),(w3),(w4) and (w5) the claim is true for $q=k+r\/$.

\vspace{1mm} Hence we conclude that $v_{k+r}x\/$ is not a radial edge of $G\/$ for any $r=0, 1, \ldots, k-3\/$.

\vspace{3mm} We now prove the sufficiency. Without loss of generality assume that the radial edges are $v_ix$ and $v_ix\in E(F_{i+1})$, $i=0,1,2,...,k-1$.\

\vspace{2mm} Consider the edge $v_1v_q$ in $F_1$. Clearly $q\notin \{2k, 2k+1,...,4k-4\}\/$ because $F_1$ is a non-crossing perfect matching. Further, $q$ is an even integer; otherwise $v_1v_q$ separates $V(G)-\{v_1,v_q\}$ into two parts each with an odd number of vertices (which is impossible by Proposition \ref{s-even} since $F_1$ is a non-crossing perfect matching). If $q< 2k-2$, then by Lemma \ref{wheel-l1}(ii),  $F_1$ contains a boundary edge which joins two vertices $v_i,   \ v_{i+1}$ in $N(x)$ (yielding a triangle $xv_iv_{i+1}x$ in G) which is impossible. Hence $q=2k-2\/$.

\vspace{2mm} By a similar recursive argument, we see that $v_jv_{2k-j-1}$ is an edge in $F_1$ for each $j=1,2,...,k-1\/$ with $v_{k-1}v_k\/$ being a boundary edge of $F_1$ (by Lemma \ref{wheel-l1}(ii)).

\vspace{2mm} By repeating the same argument to $F_2,...,F_k$ successively, we see that, for each $i=2,...,k$, $v_{i+j-1}v_{2k+i-j-2}$ is an edge of $F_i$ for each $j=1,2,...,k-1\/$ with $v_{k+i-2}v_{k+i-1}$ being a boundary edge of $F_i$ (by Lemma \ref{wheel-l1}(ii)).

\vspace{2mm}
It remains to show that each $F_i$ has just one more boundary edge. Consider $F_k$ first. Since $v_{k-2}v_p$ is an edge in $F_k$ for some $p\in \{3k-2, 3k-1,...,4k-4\}$, we apply similar argument as before (which was done to the case $F_1$ and $v_1v_q$) to conclude that $p=3k-2$. Continue with the same argument, it follows that $v_{k-1-j}v_{3k-3+j}$ is an edge of $F_k$ for each $j=1,2,...,k-1$ with $v_0v_{4k-4}$ being a boundary edge of $F_k$ (by Lemma \ref{wheel-l1}(ii)).

Now repeat the same argument to $F_{k-1},...,F_1\/$ successively, we see that, for each $i=1,2,...,k-1$, $v_{k-1-i-j}v_{3k-3-i+j}$ is an edge of $F_{k-i}$ for each $j=1,2,...,k-1\/$ with $v_{-i}v_{4k-4-i}$ being a boundary edge of $F_{k-i}$ (by Lemma \ref{wheel-l1}(ii)).

\vspace{1mm} This completes the proof. \qed

\vspace{2mm}
\begin{note}
When  $n \geq 4\/$ is even, Theorem \ref{w-oddn2} remains true. Here we may assume that $F_1 = \{ v_0x, v_jv_{n-j+1}, v_{n+j}v_{-j} \ | \ j=1, 2, \ldots, n/2 \}\/$. Note that  we may also assume that $F_1 = \{ v_0x, v_jv_{n-j-1}, v_{n-2+j}v_{-j} \ | \ j=1, 2, \ldots, n/2 \}\/$ (as the  two graphs constructed are isomorphic). The proof is similar to the case when $n\/$ is odd (with suitable modification).
\end{note}

\section{Points in $R\/$-Position}  \label{rposition}

\begin{define}
Let $S\/$ be a set of $2n\/$ points in general position (that is, no three points are collinear). We say that the points in $S\/$ are in {\em $R\/$-position\/} if there is a set $\cal L\/$ of pairwise non-parallel lines with exactly one point of $S\/$ in each open unbounded region formed by $\cal L\/$.
\end{define}

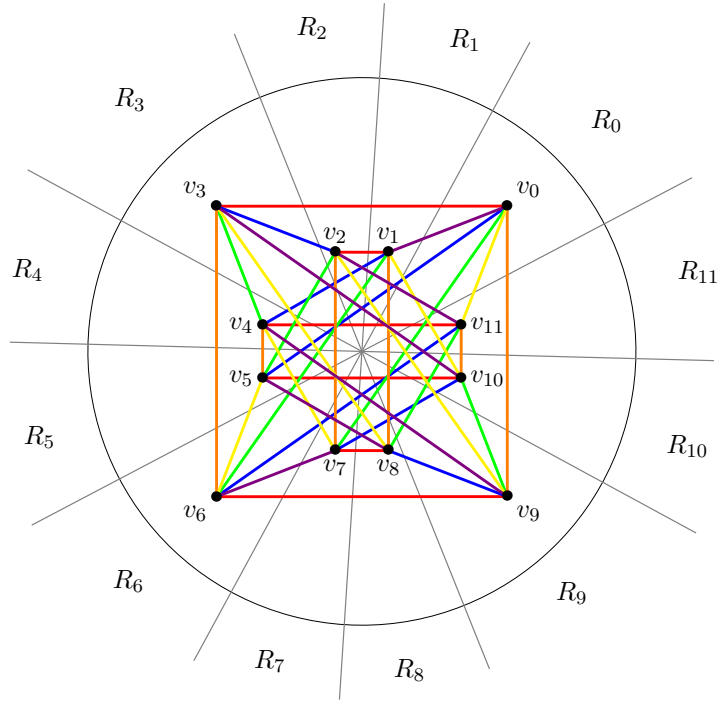
\begin{figure}[htb]
\centering
\resizebox{10cm}{!}{
\begin{tikzpicture}[rotate=75,.style={draw}]
\coordinate (center) at (0,0);
  \def\radius{4cm}
   \draw (center) circle[radius=\radius];
 %  \foreach \x in {0,27.6923,...,360} {
  %          \filldraw[] (\x:3.5cm) circle(1pt);
%             \filldraw[] (0,0) circle(1pt);
  %             }
        \draw [line width=0.5,gray](-1.2,5) -- (1.2,-5);\draw [line width=0.5,gray](5,1) -- (-5,-1);\draw [line width=0.5,gray](4,3) -- (-4,-3);\draw [line width=0.5,gray](1.3,5.4) -- (-1.3,-5.4);\draw [line width=0.5,gray](5,-1.2) -- (-5,1.2);\draw [line width=0.5,gray](3.7,-4) -- (-3.7,4);
         %points
       \coordinate (v1) at (1.5,0);\node at (v1) {\textbullet};\node[above] at (v1) {${v_1}$};
       \coordinate (v2) at (1.3,0.75);\node at (v2) {\textbullet};\node[above] at (v2) {${v_2}$};
       \coordinate (v3) at (1.5,2.6);\node at (v3) {\textbullet};\node[above left] at (v3) {${v_3}$};
       \coordinate (v4) at (0,1.5);\node at (v4) {\textbullet};\node[left] at (v4) {${v_4}$};
       \coordinate (v5) at (-0.75,1.3);\node at (v5) {\textbullet};\node[left] at (v5) {${v_5}$};
       \coordinate (v6) at (-2.6,1.5);\node at (v6) {\textbullet};\node[below left] at (v6) {${v_6}$};
       \coordinate (v7) at (-1.5,0);\node at (v7) {\textbullet};\node[below] at (v7) {${v_7}$};
       \coordinate (v8) at (-1.3,-0.75);\node at (v8) {\textbullet};\node[below] at (v8) {${v_8}$};
       \coordinate (v9) at (-1.5,-2.6);\node at (v9) {\textbullet};\node[below right] at (v9) {${v_9}$};
       \coordinate (v10) at (0,-1.5);\node at (v10) {\textbullet};\node[right] at (v10) {${v_{10}}$};
       \coordinate (v11) at (0.75,-1.3);\node at (v11) {\textbullet};\node[right] at (v11) {${v_{11}}$};
       \coordinate (v0) at (2.6,-1.5);\node at (v0) {\textbullet};\node[above right] at (v0) {${v_{0}}$};
       %%%%%%%%%%%%%%%%%%%%%%%%%%%%%%%%%%%%%%%%%%%%%%%%%%%%%%%%%%%%%%%%%%%%%%%%%%%%%%%%%%%%%%%%%%%%%%%%%%%%
        \coordinate (R0) at (3.8,-2.3);\node at (R0) {};\node[above right] at (R0) {${R_{0}}$};
        \coordinate (R1) at (4.4,0);\node at (R1) {};\node[above right] at (R1) {${R_{1}}$};
         \coordinate (R2) at (4,2.2);\node at (R2) {};\node[above right] at (R2) {${R_{2}}$};
         \coordinate (R3) at (2.3,4.5);\node at (R3) {};\node[above right] at (R3) {${R_{3}}$};
         \coordinate (R4) at (-0.5,5.3);\node at (R4) {};\node[above right] at (R4) {${R_{4}}$};
         \coordinate (R5) at (-2.8,4.5);\node at (R5) {};\node[above right] at (R5) {${R_{5}}$};
         \coordinate (R6) at(-4.5,2.7);\node at (R6) {};\node[above right] at (R6) {${R_{6}}$};
         \coordinate (R7) at(-5.1,0.4);\node at (R7) {};\node[above right] at (R7) {${R_{7}}$};
        \coordinate (R8) at(-4.7,-1.6);\node at (R8) {};\node[above right] at (R8) {${R_{8}}$};
        \coordinate (R9) at(-3,-3.6);\node at (R9) {};\node[above right] at (R9) {${R_{9}}$};
        \coordinate (R10) at(-0.5,-4.6);\node at (R10) {};\node[above right] at (R10) {${R_{10}}$};
        \coordinate (R11) at(2,-4.1);\node at (R11) {};\node[above right] at (R11) {${R_{11}}$};

       %%%%%%%%%%%%%%%%%%%%%%%%%%%%%%%%%%%%%%%%%%%%%%%%%%%%%%%%%%%%%%%%%%%%%%%%%%%%%%%%%%%%%%%%%%%%%%%%%%%%
         %lines
       \draw [line width=1.2,red](v1) -- (v2);\draw [line width=1.2,red](v0) -- (v3);\draw [line width=1.2,red](v11) -- (v4);
       \draw [line width=1.2,red](v10) -- (v5);\draw [line width=1.2,red](v9) -- (v6); \draw [line width=1.2,red](v8) -- (v7);
       \draw [line width=1.2,blue](v2) -- (v3); \draw [line width=1.2,blue](v1) -- (v4); \draw [line width=1.2,blue](v0) -- (v5);
       \draw [line width=1.2,blue](v11) -- (v6);\draw [line width=1.2,blue](v10) -- (v7);\draw [line width=1.2,blue](v9) -- (v8);
       \draw [line width=1.2,green](v3) -- (v4); \draw [line width=1.2,green](v2) -- (v5); \draw [line width=1.2,green](v1) -- (v6);
       \draw [line width=1.2,green](v0) -- (v7); \draw [line width=1.2,green](v11) -- (v8); \draw [line width=1.2,green](v10) -- (v9);
       \draw [line width=1.2,orange](v4) -- (v5); \draw [line width=1.2,orange](v3) -- (v6); \draw [line width=1.2,orange](v2) -- (v7);
       \draw [line width=1.2,orange](v1) -- (v8); \draw [line width=1.2,orange](v0) -- (v9); \draw [line width=1.2,orange](v11) -- (v10);
       \draw [line width=1.2,yellow](v5) -- (v6); \draw [line width=1.2,yellow](v4) -- (v7); \draw [line width=1.2,yellow](v3) -- (v8);
       \draw [line width=1.2,yellow](v2) -- (v9); \draw [line width=1.2,yellow](v1) -- (v10); \draw [line width=1.2,yellow](v0) -- (v11);
       \draw [line width=1.2,violet](v6) -- (v7); \draw [line width=1.2,violet](v5) -- (v8); \draw [line width=1.2,violet](v4) -- (v9);
       \draw [line width=1.2,violet](v3) -- (v10); \draw [line width=1.2,violet](v2) -- (v11); \draw [line width=1.2,violet](v1) -- (v0);

 \node at (v1) {\textbullet};
      \node at (v2) {\textbullet};
      \node at (v3) {\textbullet};
      \node at (v4) {\textbullet};
      \node at (v5) {\textbullet};
      \node at (v6) {\textbullet};
      \node at (v7) {\textbullet};
      \node at (v8) {\textbullet};
      \node at (v9) {\textbullet};
      \node at (v10) {\textbullet};
      \node at (v11) {\textbullet};
      \node at (v0) {\textbullet};

\end{tikzpicture}
}
\caption{ Triangle-free geometric graph in $R\/$-position with $6\/$ non-crossing perfect matchings } \label{sufficint}

\end{figure}

In \cite{bhrw:refer}, it was noted that if the points in $S\/$ are in $R\/$-position, then the unbounded regions and the  points in $S\/$ can be labeled as $R_0, R_1, \ldots, R_{2n-1}\/$ and    $v_0, v_1, \ldots, v_{2n-1}\/$ in anti-clockwise direction respectively with $R_i\/$ containing $v_i\/$, $i=0, 1, \ldots, 2n-1\/$ (see Figure \ref{sufficint}).  Also, the authors showed that if the vertex set of the complete geometric graph $K_{2n}\/$ are in the $R\/$-position, then the edge-set of $K_{2n}\/$ can be partitioned into $n\/$ plane spanning double stars (which are pairwise graph-isomorphic).

\vspace{1mm} When the points of $S\/$ are in $R\/$-position, we have the following condition for the existence of $n\/$ non-crossing perfect matchings whose union is a maximal triangle-free graph. An example of such graph with $n=6\/$ is depicted in Figure \ref{sufficint}.

\begin{result}  \label{r-post}
Let $S\/$ be a set of $2n\/$ points $v_0, v_1, \ldots, v_{2n-1}\/$  in $R\/$-position where $n \geq 2\/$. Suppose $i \in \{0, 1, \ldots, 2n-1\}\/$ and $j \leq i+n \ (mod \ 2n)\/$ are of different parity, and $v_i\/$ and $v_j\/$ are not separated by the line passing through $v_{i-1}, v_{j+1}\/$.   Then there exist $n\/$ non-crossing perfect matchings whose union is a maximal triangle-free geometric graph.
\end{result}

\vspace{1mm}  \noindent
{\bf Proof:}
For each $i = 0,1,\dots , n-1\/$, let $F_i\/$ denote the set of all edges $v_lv_k$ where $l+k\equiv 2i+1 \ (mod \ 2n)\/$. That is, any two edges $v_r v_s\/$ and $v_lv_k\/$ are in the same set $F_i\/$ if and only if $r+s\equiv l+k \ (mod \ 2n)\/$, see Figure \ref{sufficint}.\

\vspace{1mm} Clearly, by definition of $F_i$, the sets $F_0, F_1, \ldots, F_{n-1}$ are perfect matchings on $S\/$. Moreover $E(F_i)\cap E(F_j)= \phi\/$ for $i\neq j\/$.\

\vspace{1mm} Assume that $v_lv_k \in F_i\/$ for some $0\leq i\leq n-1$. Then by definition of $F_i$,   $l$ and $k$ are of different parity. It is clear that $v_{l-1}v_{k+1} \in F_i\/$ since $(l-1)+(k+1)\equiv l+k\/$.

\vspace{1mm} Note that the  two edges $v_lv_k\/$ and $v_{l-1}v_{k+1}\/$ are non-crossing, otherwise, $v_l\/$ and $v_k\/$ are separated by the line passing through $v_{l-1}, v_{k+1}\/$, a contradiction. Similarly, $v_{l-1}v_{k+1}\/$ and $v_{l-2}v_{k+2}\/$ are non-crossing edges.\

\vspace{1mm} Suppose on the contrary  that $v_lv_k\/$ and $v_{l-2}v_{k+2}\/$ are crossing. That is,  $v_l\/$ and $v_k\/$ are separated by the line passing through $v_{l-2}\/$ and $v_{k+2}\/$. We consider the location of $v_{k-1}\/$.

\begin{figure}[htb]
\resizebox{15cm}{!}{
\begin{minipage}{0.75\textwidth}
\begin{tikzpicture}
%\coordinate (center) at (0,0);
%  \def\radius{2.5cm}
%   \draw (center) circle[radius=\radius];
%   \foreach \x in {0,45,...,360} {
%             \filldraw[] (\x:2.5cm) circle(1pt);
%             \filldraw[] (0,0) circle(1pt);
%               }
  %               \draw [line width=1,green](-1.8,1.8) -- (0,0); \draw [line width=1,green](1.8,-1.8) -- (0,0);

           %points
       \coordinate (vk) at (-0.5,-0.5);\filldraw[black] (vk) circle(2pt);\node[above] at (vk) {${v_{k}}$};
       \coordinate (vl-1) at (3,1.8);\filldraw[black] (vl-1) circle(2pt);\node[above right] at (vl-1) {${v_{l-1}}$};
       \coordinate (vl+1) at (0.8,1.3);\filldraw[black] (vl+1) circle(2pt);\node[above right] at (vl+1) {${v_{l+1}}$};
       \coordinate (vl-2) at (2.5,0);\filldraw[black] (vl-2) circle(2pt);\node[right] at (vl-2) {${v_{l-2}}$};
       \coordinate (vl) at (1.4,0.5);\filldraw[black] (vl) circle(2pt);\node[above] at (vl) {${v_l}$};
       \coordinate (vk-1) at (-0.5,1.3);\filldraw[black] (vk-1) circle(2pt);\node[above] at (vk-1) {${v_{k-1}}$};
       \coordinate (vk+1) at (-2.8,1.8);\filldraw[black] (vk+1) circle(2pt);\node[above left] at (vk+1) {${v_{k+1}}$};
       \coordinate (vk+2) at (-2.5,0);\filldraw[black] (vk+2) circle(2pt);\node[left] at (vk+2) {${v_{k+2}}$};
%%%%%%%%%%%%%%%%%%%%%%%%%%%%%%%%%%%%%%%%%%%%%%%%%%%%%%%%%%%%%%%%%%%%%%%%%%%%%%%%%%%%%%%%%%%%%%%%%%%%%%%%%%%%%%%%%%%%%%%
\coordinate (vll) at (2.5,3.3);\coordinate (vkk) at (-4.1,-1.1);
%%%%%%%%%%%%%%%%%%%%%%%%%%%%%%%%%%%%%%%%%%%%%%%%%%%%%%%%%%%%%%%%%%%%%%%%%%%%%%%%%%%%%%%%%%%%%%%%%%%%%%%%%%%%%%%%%%%%%%
\coordinate (x) at (0,-1.7);
\coordinate (v0) at (4.7,3);\node[below right] at (v0) {${R_{l-2}}$};
\coordinate (v1) at (3.8,3.4);\node[below right] at (v1) {${R_{l-1}}$};
\coordinate (v2) at (2.5,3.9);\node[below right] at (2.7,3.9) {${R_{l}}$};
\coordinate (v3) at (1.1,3.9);\node[below right] at (1.1,4.1) {${R_{l+1}}$};
\coordinate (v4) at (-0.3,4);
\coordinate (v5) at (-1.4,3.9);\node[below right] at (-1.4,4.1) {${R_{k-1}}$};
\coordinate (v6) at (-3,3.5);\node[below right] at  (-2.7,4){${R_{k}}$};
\coordinate (v7) at (-4.2,2.6);\node[below right] at  (-4,3.4) {${R_{k+1}}$};
\coordinate (v8) at (-4.2,2.6);\node[below left] at  (-4.2,2.6) {${R_{k+2}}$};
   %lines
       \draw [line width=1,black](vl-2) -- (vk+2); \draw [line width=1,black](vl-1) -- (vk+1); \draw [line width=1,black](vl+1) -- (vk-1); \draw [line width=1,black](vl) -- (vk);
%%%%%%%%%%%%%%%%%%%%%%%%%%%%%%%%%%%%%%%%%%%%%%%%%%%%%%%
\draw [dashed,gray](vll) -- (vkk);\draw [line width=0.5,gray](x) -- (v0);\draw [line width=0.5,gray](x) -- (v1);\draw [line width=0.5,gray](x) -- (v2);\draw [line width=0.5,gray](x) -- (v3);\draw [line width=0.5,gray](x) -- (v4);\draw [line width=0.5,gray](x) -- (v5);\draw [line width=0.5,gray](x) -- (v6);\draw [line width=0.5,gray](x) -- (v7);

   %    \node at (v1) {\textbullet};
  %   \node at (x) {\textbullet}; \filldraw[green] (v0) circle(2pt);\filldraw[green] (v1) circle(2pt);\filldraw[green] (v2) circle(2pt);

\end{tikzpicture}
\centering

Case (1)

\end{minipage}
\begin{minipage}{0.8\textwidth}
\begin{tikzpicture}
%\coordinate (center) at (0,0);
%  \def\radius{2.5cm}
%   \draw (center) circle[radius=\radius];
%   \foreach \x in {0,45,...,360} {
%             \filldraw[] (\x:2.5cm) circle(1pt);
%             \filldraw[] (0,0) circle(1pt);
%               }
  %               \draw [line width=1,green](-1.8,1.8) -- (0,0); \draw [line width=1,green](1.8,-1.8) -- (0,0);

           %points
       \coordinate (vk) at (-0.5,-0.5);\filldraw[black] (vk) circle(2pt);\node[above] at (vk) {${v_{k}}$};
       \coordinate (vl-1) at (3,1.8);\filldraw[black] (vl-1) circle(2pt);\node[above right] at (vl-1) {${v_{l-1}}$};
       \coordinate (vl+1) at (0.3,-0.7);\filldraw[black] (vl+1) circle(2pt);\node[above right] at (vl+1) {${v_{l+1}}$};
       \coordinate (vl-2) at (2.5,0.2);\filldraw[black] (vl-2) circle(2pt);\node[right] at (vl-2) {${v_{l-2}}$};
       \coordinate (vl) at (1.4,0.9);\filldraw[black] (vl) circle(2pt);\node[above right] at (vl) {${v_l}$};
       \coordinate (vk-1) at (-0.16,-0.7);\filldraw[black] (vk-1) circle(2pt);\node[below] at (vk-1) {${v_{k-1}}$};
       \coordinate (vk+1) at (-2.8,1.8);\filldraw[black] (vk+1) circle(2pt);\node[above left] at (vk+1) {${v_{k+1}}$};
       \coordinate (vk+2) at (-2.5,0.2);\filldraw[black] (vk+2) circle(2pt);\node[left] at (vk+2) {${v_{k+2}}$};
%%%%%%%%%%%%%%%%%%%%%%%%%%%%%%%%%%%%%%%%%%%%%%%%%%%%%%%%%%%%%%%%%%%%%%%%%%%%%%%%%%%%%%%%%%%%%%%%%%%%%%%%%%%%%%%%%%%%%%%
\coordinate (vll) at (4.5,1.5);\coordinate (vkk) at (-4.1,-0.14);
%%%%%%%%%%%%%%%%%%%%%%%%%%%%%%%%%%%%%%%%%%%%%%%%%%%%%%%%%%%%%%%%%%%%%%%%%%%%%%%%%%%%%%%%%%%%%%%%%%%%%%%%%%%%%%%%%%%%%%
\coordinate (x) at (0,-1.7);
\coordinate (v0) at (4.7,3);\node[below right] at (v0) {${R_{l-2}}$};
\coordinate (v1) at (3.8,3.4);\node[below right] at (v1) {${R_{l-1}}$};
\coordinate (v2) at (2.5,3.9);\node[below right] at (2.7,3.9) {${R_{l}}$};
\coordinate (v3) at (1.1,3.9);\node[below right] at (1.1,4.1) {${R_{l+1}}$};
\coordinate (v4) at (-0.3,4);
\coordinate (v5) at (-1.4,3.9);\node[below right] at (-1.4,4.1) {${R_{k-1}}$};
\coordinate (v6) at (-3,3.5);\node[below right] at  (-2.7,4){${R_{k}}$};
\coordinate (v7) at (-4.2,2.6);\node[below right] at  (-4,3.4) {${R_{k+1}}$};
\coordinate (v8) at (-4.2,2.6);\node[below left] at  (-4.2,2.6) {${R_{k+2}}$};
   %lines
       \draw [line width=1,black](vl-2) -- (vk+2); \draw [line width=1,black](vl-1) -- (vk+1); \draw [line width=1,black](vl+1) -- (vk-1); \draw [line width=1,black](vl) -- (vk);
%%%%%%%%%%%%%%%%%%%%%%%%%%%%%%%%%%%%%%%%%%%%%%%%%%%%%%%
\draw [dashed,gray](vll) -- (vkk);\draw [line width=0.5,gray](x) -- (v0);\draw [line width=0.5,gray](x) -- (v1);\draw [line width=0.5,gray](x) -- (v2);\draw [line width=0.5,gray](x) -- (v3);\draw [line width=0.5,gray](x) -- (v4);\draw [line width=0.5,gray](x) -- (v5);\draw [line width=0.5,gray](x) -- (v6);\draw [line width=0.5,gray](x) -- (v7);

   %    \node at (v1) {\textbullet};
  %   \node at (x) {\textbullet}; \filldraw[green] (v0) circle(2pt);\filldraw[green] (v1) circle(2pt);\filldraw[green] (v2) circle(2pt);

\end{tikzpicture}
\centering

Case (2)

\end{minipage}
}
\caption{Points in $R\/$-position. }  \label{case2}
\end{figure}
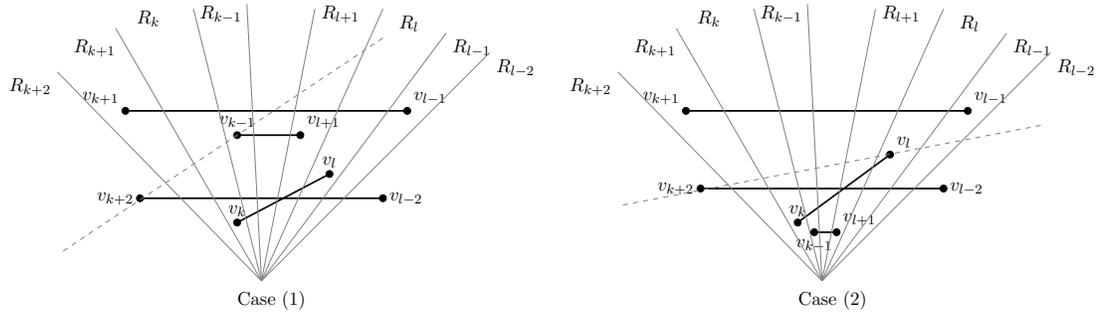

\vspace{1mm} (a) Suppose $v_{k-1}\/$ is in the region $R_{k-1}\/$ bounded by the line passing through $v_{k+2}\/$ and $v_{l-2}\/$ and the line passing through $v_{k+1}\/$ and $v_{l-1}\/$ (see Case (1) of  Figure \ref{case2}).

\vspace{2mm}

\vspace{1mm} In this case, $v_k, v_{k+1}\/$ are separated by the line passing through $v_{k-1}\/$ and $v_{k+2}\/$, a contradiction.

\vspace{1mm} (b)  Suppose $v_{k-1}\/$ is not in the region $R_{k-1}\/$ bounded by the line passing through $v_{k+2}\/$ and $v_{l-2}\/$ and the line passing through $v_{k+1}\/$ and $v_{l-1}\/$ (see Case (2) of  Figure \ref{case2}).

\vspace{1mm} (i) $v_{k-1}, v_k\/$ is not separated by the line passing through $v_{l-2}\/$ and $v_{k+2}\/$. Here $v_{l+2}\/$ and $v_k\/$ are separated by the line passing through $v_{l+1}\/$ and $v_{k+1}\/$, a contradiction.

\vspace{1mm} (ii) $v_{k-1}, v_k\/$ is separated by the line passing through $v_{l-2}\/$ and $v_{k+2}\/$. Here $v_{k-1}\/$ and $v_k\/$ are separated by the line passing through $v_{k-2}\/$ and $v_{k+1}\/$, a contradiction.

\vspace{1mm} We may color the vertices $v_{2i} \/$ with one color and the vertices $v_{2i+1} \/$ with another color, $i = 0,1, \ldots, n-1\/$. It is easy to see that $\bigcup_{i=0} ^{n-1} F_i\/$ is a triangle-free graph having $2n\/$ vertices and $n^{2}\/$ edges (since no edge joins two vertices of the same color). By Turan's theorem, $\bigcup_{i=0} ^{n-1} F_i\/$ is a maximal triangle-free graph. \qed

\begin{note}
Note that, when the  $2n\/$ points of $S\/$ are in convex position, then $S\/$ clearly satisfies the conditions of Theorem \ref{r-post}. Hence there exist $n\/$ edge-disjoint non-crossing perfect matchings whose union is a maximal triangle-free graph (Theorem \ref{c2n-n}). However, when the points of $S\/$ are in regular wheel configuration, then $S\/$ does not satisfy the conditions of Theorem \ref{r-post}.
\end{note}

\section{Points in general position}   \label{gposition}

Let $S\/$ be a set of $2n\/$ points in  plane where no three points  are collinear. In \cite{lms:refer}, it was mentioned that if $S\/$ is a disjoint union of $S_1, S_2\/$, then the ham-sandwich theorem guarantees the existence of a line that bisects $S_1\/$ and $S_2\/$.     In this section, for $n =2^k +h\/$ where $0 \leq h < 2^k\/$, we shall show that there are  at least $k\/$  edge-disjoint non-crossing perfect matchings such that the union is a triangle-free geometric graph (Theorem \ref{general}).

 \vspace{1mm} For this purpose, we recall a result of \cite{ap:refer} which states that if $V$ is a set of $n\/$ points in general position in the plane with $n = n_1 + n_2 + \cdots + n_{2m}$, where the $n_i's$ are nonnegative integers, then there exist $m$ lines  and a partition of $V$ into $2m$ pairwise disjoint subsets $V_1,V_2,...,V_{2m}$, such that $|V_i|=n_i$, and every two distinct subsets $V_i,V_j$ are separated by at least one of the $m$ lines. A special case of this result is the following.

\begin{deduce} \cite{ap:refer}
Let $S\/$ be a set of $n\/$ points in general position in the plane with $n = n_1 + n_2\/$. Then there exists a line that separates $S\/$ into two disjoint sets $S_1\/$ and $S_2\/$ such that $|S_1| =n_1\/$ and $|S_2| = n_2\/$.
\end{deduce}

In the case that $|S_1| = |S_2|\/$, there is an algorithm due to Attalah (see \cite{ata:refer}) that matches  all points in $S_1\/$  to all points in $S_2\/$ resulting in a non-crossing perfect matching in $S\/$.  In the case that $|S_1| = |S_2| +2\/$, we shall modify this algorithm to construct a non-crossing perfect matching $F\/$ in $S\/$. The modified algorithm is given below where $|S| = 2n\/$ and $n =2^k +h\/$.

\vspace{3mm}
 Algorithm $(A):\/$
\begin{enumerate}
\item Find a line $l\/$ that separates $S\/$ into $S_1\/$ and $S_2\/$ such that   either $|S_1| = |S_2|\/$ or $|S_1| = |S_2| +2\/$.

\item Find a line $l^{\bot}\/$ such that $l^{\bot}\/$ is perpendicular to $l\/$ and all points in $S\/$ are on one side of  $l^{\bot}\/$.

\item  Find $CH(S_i)\/$, the convex hull of $S_i\/$, $i=1, 2\/$ and let $u_i \in CH(S_i)\/$, $i=1, 2\/$ be such that all the points in $S_1 \cup S_2 -\{u_1, u_2\}\/$ are between $l^{\bot}\/$ and  the line joining $u_1\/$ and $u_2\/$.  Let $u_1u_2\/$ be an edge in $F\/$.

\item Repeat Step $3\/$ with $S_i - u_i\/$ taking the place of $S_i\/$, $i=1,2\/$ until $S_2 = \emptyset\/$.

\item  If (i) $S_1 =\emptyset\/$, then stop; otherwise (ii)  $S_1 = \{v_1, v_2\}\/$ then let $v_1v_2\/$ be an edge in $F\/$ and stop.

\end{enumerate}

\vspace{1mm} The edge  $v_1v_2\/$ in Step 5(ii) is termed a {\em stone\/} and shall be  denoted by $st(v_1, v_2)\/$.

\vspace{1mm}
\begin{support} \label{st}
Let $T\/$ be a set of $m\/$ points in the general position where $m \geq 3\/$. Suppose there is a line separating a given set  $\{u, v\}\/$  from $T\/$. Then there is a line that separates some subset $T_1\/$ from $T\/$ with $\{u, v\} \subseteq T_1\/$ and $2 \leq |T_1| \leq m-1\/$.
\end{support}

\vspace{1mm}  \noindent
{\bf Proof:} The lemma is trivially true if $m=3\/$ with $T_1 = \{u, v\}\/$. So assume that $m \geq 4\/$.

\vspace{1mm}
Let $w_1 \in T - \{u, v\}\/$ be such that all points in $T - \{u, v, w_1\}\/$ are on one side of the line $l_1\/$ joining $w_1\/$ and $z_1\/$ for some $z_1 \in \{u, v\}\/$ and let $T_1 = \{u, v, w_1\}\/$. Let $L_1\/$ be a line parallel to $l_1\/$ such that all points in $T - \{u, v, w_1\}\/$ are on one side of $L_1\/$.

\vspace{1mm} If $|T_1| = m-1\/$, then the proof is complete.  Otherwise  repeat the argument with $w_2 \in T - \{u, v, w_1\}\/$ and $z_2 \in  \{u, v, w_1\}\/$ so that all points in $T - \{u, v, w_1, w_2\}\/$ are on one side of the line $l_2\/$ joining $w_2z_2\/$ and let $T_1 = \{ u, v, w_1, w_2\}\/$ with the line $L_2\/$ similarly defined. By repeating the argument where necessary, we reach the conclusion of the lemma. \qed

\vspace{2mm} We shall now apply Algorithm $(A)\/$  to prove the next result. We wish to emphasize that in Step $1\/$, we find the line that separates $S\/$ into $S_1\/$ and $S_2\/$ so that $|S_i|\/$ is even for each $i=1, 2\/$ in  each iteration.

\vspace{1mm}
\begin{result} \label{general}
Let $S$ be a set of $2n$ points in general position in the plane where $n=2^k+h$, with $0\leq h <2^k$. Then
there exist at least $k$  edge-disjoint non-crossing perfect matchings $ F_1, F_2, ... , F_k$ such that $\bigcup_{i=1} ^k F_i\/$ is a triangle-free geometric graph .
\end{result}

\vspace{1mm}  \noindent
{\bf Proof:} First we apply Algorithm $(A)\/$  above to obtain the first non-crossing perfect matching $F_1\/$. In so doing, the set $S\/$ has been split into $S_1\/$ and $S_2\/$ and $|S_i|\/$ is even, $i=1, 2\/$.

\vspace{1mm} If $S_i\/$ has no stone, then we apply Algorithm $(A)\/$ to split $S_i\/$ into $S_{i,1}\/$ and $S_{i,2}\/$ and  obtain an non-crossing perfect matching $F(i)\/$ on $S_i\/$ and let $F_2 = F(1) \cup F(2)\/$.

\vspace{1mm} If $S_i\/$ has a stone $st(v_1, v_2)\/$, then Algorithm $(A)\/$ ensures that there is a line separating $\{v_1, v_2\}\/$ from $S_i\/$. Hence by Lemma \ref{st}, there is line that splits $S_i\/$ into $S_{i,1}\/$ and $S_{i, 2}\/$ with $\{v_1, v_2\} \subseteq S_{i,2}\/$,  $|S_{i,j}|\/$ is even for $j =1, 2\/$ and either $|S_{i, 2}| = |S_{i, 1}|\/$ or  $|S_{i, 2}| = |S_{i, 1}|-2\/$. Let $F(i) \/$ be a non-crossing perfect matching in $S_i\/$ and let $F_2 = F(1) \cup F(2)\/$.

\vspace{1mm}
To obtain the next non-crossing perfect matching, we repeat the above operations by applying Algorithm $(A)\/$ to $S_{i,j}\/$ for each $i=1, 2\/$ with $j = 1, 2\/$ to obtain a non-crossing perfect matching $F(i)\/$ for $S_i = S_{i, 1} \cup S_{i,2}\/$ and let $F_3 = F(1) \cup F(2)\/$.

\vspace{1mm}
Continue with the above operations until we obtain $k\/$ sets  of non-crossing perfect matchings $F_1, F_2, \ldots, F_k\/$. It is clear that $E(F_i) \cap E(F_j) = \emptyset\/$ for $i \neq j\/$.

\vspace{1mm} Next we show that $G=\bigcup_{i=1} ^k F_i\/$ is a triangle-free geometric graph.

\vspace{1mm} Assume on the contrary  that G has a triangle $\Delta=v_rv_sv_tv_r$. Then at most one of the edges in $\Delta\/$ is a stone (by the above construction).

\vspace{1mm} (i) Suppose $st(v_r, v_s)\/$ is a stone. % Then $v_rv_s\/$ belongs to some non-crossing perfect matching $F_i\/$.
Since $v_rv_t\/$ is an edge in $\Delta\/$, $v_rv_t\/$ belongs to some non-crossing perfect matching $F_j\/$.  By the above construction, there is a line $l'\/$ which separates $v_r\/$ and $v_t\/$ with $v_r, v_s\/$ on the the same side of $l'\/$. But this implies that $v_t\/$ cannot be adjacent to any vertex  which lies in the same side as $v_r, v_s\/$ (with respect to $l'\/$), a contradiction.

\vspace{1mm}
(ii) Suppose no edge in $\Delta\/$ is a stone. Again $v_rv_t\/$ belongs to some non-crossing perfect matching $F_j\/$.   By the above construction, there is a line $l'\/$ which separates $v_r\/$ and $v_t\/$. We can assume without loss o generality that $v_r,  v_s\/$ are on the same side of $l'\/$. But again we reach the same contradiction as in Case (i).

\vspace{1mm} This completes the proof.   $\qed$

%%%%%%%%%%%%%%%%%%%%%%%%%%%%%%%%%%%%%%%%%%%%%%%%%%%%%%%%%%%%%%%%%%%%%%%%%%%%%%%%%%%%%%%%%%%%%%%%%%%%%%%%%%%%%%%%%%%%%%%%%%%%%%%
%%%%%%%%%%%%%%%%%%%%%%%%%%%%%%%%%%%%%%%%%%%%%%%%%%%%%%%%%%%%%%%%%%%%%%%%%%%%%%%%%%%%%%%%%%%%%%%%%%%%%%%%%%%%%%%%%%%%%%%%%%%%%%

%\vspace{1mm} Suppose that there is $st(u, v) $ in $\Delta=v_rv_sv_tv_r$. Without loss of generality assume that $st(u, v) = v_r v_s$. Now $v_rv_s\in F_{i_1}$ implies that $v_r\in P_{i_1,j}$ and $v_s\in P_{i_1,j}$. But $v_sv_t\in F_{i_2}$  implies that $v_s\in P_{i_2,j^*}$ and $v_t\in P_{i_2,j^* +1}$ , and $v_rv_t\in F_{i_3}$ implies that $v_r$ and $v_t$ must be in the same part in the previous iteration $ i_2 $ so $v_r\in P_{i_2,j^* +1}$ but $v_r v_s = st(v_r, v_s)$ and that implies $v_s\in P_{i_2,j^* +1}$  and this contradiction (since $P_{i_2,j^*}\cap P_{i_2,j^* +1}=\phi$). $\qed$\\

\section{Plane Triangle-free Geometric Graphs}

In \cite{bbms:refer} (Theorem 8), the authors prove that for a set of $2n\/$ points $S\/$ in the general position, where $n \geq 2\/$,  there exist at least $2\/$  and at most $5\/$ edge-disjoint non-crossing perfect matchings that can be packed into  a complete geometric  graph $K_{|S|}\/$ on the set $S\/$. Moreover these bounds are tight.

\vspace{1mm}
For the case of triangle-free geometric graphs in general position, we have the following. The proof follows easily from the fact that every triangle-free simple planar graph has a vertex of degree at most $3\/$.

\vspace{1mm}
\begin{propo} \label{3pm}
 Suppose $S\/$  is a set of $2n\/$ points in the plane where $n \geq 4\/$.
 If $S\/$ is in general position, then at most $3\/$  edge-disjoint non-crossing perfect matchings can be packed into $K_{|S|}\/$ such that the union of these perfect matchings is a triangle-free plane geometric graph. On the the hand, if $S\/$ is in convex position, then at most $2\/$  such  perfect matchings can be packed into $K_{|S|}\/$ giving rise to a triangle-free plane geometric graph.
\end{propo}

\vspace{1mm}
For every natural number $n \geq 4\/$, the prism (which is the Cartesian product of an $n\/$-cycle with a path on two vertices) with $2n\/$ vertices is a triangle-free planar graph. This shows that the bound provided in Proposition \ref{3pm} is tight.

\vspace{1mm}
We conclude this section with the following result.

\begin{propo}
For a set $S\/$ of $2n \/$ points in general position in the plane, where $n \geq 4\/$, at least $2\/$ and at most $3\/$ edge-disjoint  non-crossing plane perfect matchings can be packed into $K_{|S|}\/$ so that the resulting graph is a triangle-free plane geometric graph. These bounds are tight.
\end{propo}

% \vspace{3mm}

\end{document}